\documentclass[12pt]{article}
\usepackage{fullpage}
\usepackage{amsfonts, amsmath, amsthm, stmaryrd}

\newtheorem{theorem}{Theorem}[subsection]
\newtheorem{lemma}[theorem]{Lemma}

\newtheorem{cor}[theorem]{Corollary}
\newtheorem{prop}[theorem]{Proposition}
\newtheorem{algo}[theorem]{Algorithm}
\theoremstyle{definition}
\newtheorem{defn}[theorem]{Definition}
\newtheorem{remark}[theorem]{Remark}

\newtheorem{example}[theorem]{Example}
\newtheorem{notation}[theorem]{Notation}
\numberwithin{equation}{theorem}

\newcommand{\AAA}{\mathbb{A}}
\newcommand{\CC}{\mathbb{C}}
\newcommand{\FF}{\mathbb{F}}
\newcommand{\GG}{\mathbb{G}}
\newcommand{\HH}{\mathbb{H}}

\newcommand{\PP}{\mathbb{P}}
\newcommand{\QQ}{\mathbb{Q}}
\newcommand{\ZZ}{\mathbb{Z}}
\newcommand{\calL}{\mathcal{L}}
\newcommand{\calO}{\mathcal{O}}
\newcommand{\gotho}{\mathfrak{o}}
\newcommand{\gothp}{\mathfrak{p}}
\newcommand{\frakP}{\mathfrak{P}}
\newcommand{\frakU}{\mathfrak{U}}
\newcommand{\frakX}{\mathfrak{X}}
\newcommand{\frakZ}{\mathfrak{Z}}

\newcommand{\del}{\partial}

\DeclareMathOperator{\an}{an}
\DeclareMathOperator{\Aut}{Aut}
\DeclareMathOperator{\Br}{Br}
\DeclareMathOperator{\corank}{corank}
\DeclareMathOperator{\crys}{crys}
\DeclareMathOperator{\divi}{div}
\DeclareMathOperator{\dR}{dR}
\DeclareMathOperator{\Frac}{Frac}
\DeclareMathOperator{\Gal}{Gal}

\DeclareMathOperator{\lcm}{lcm}
\DeclareMathOperator{\NS}{NS}
\DeclareMathOperator{\ord}{ord}
\DeclareMathOperator{\Pic}{Pic}

\DeclareMathOperator{\rank}{rank}
\DeclareMathOperator{\Res}{Res}
\DeclareMathOperator{\rig}{rig}
\DeclareMathOperator{\sep}{sep}
\DeclareMathOperator{\Spec}{Spec}
\DeclareMathOperator{\tors}{tors}
\DeclareMathOperator{\trace}{trace}

\begin{document}

\title{Bounding Picard numbers of surfaces using $p$-adic cohomology}
\author{Timothy G. Abbott, Kiran S. Kedlaya, and David Roe}
\date{January 18, 2007}

\maketitle

\begin{abstract}
Motivated by an application to LDPC (low density parity check)
algebraic geometry codes described by
Voloch and Zarzar,
we describe a computational procedure for establishing an upper bound on
the arithmetic or geometric Picard number
of a smooth projective surface over a finite field,
by computing the Frobenius action on $p$-adic cohomology to a small
degree of $p$-adic accuracy. 
We have implemented this procedure in \textsc{Magma}; using this 
implementation, we exhibit several examples, such as smooth quartics
over $\FF_2$ and $\FF_3$ with arithmetic Picard number 1, and a smooth quintic over
$\FF_2$ with geometric Picard number 1. We also produce some examples of
smooth quartics with geometric Picard number 2, which by a construction
of van Luijk also have trivial geometric automorphism group.
\end{abstract}

\section*{Introduction}

Much recent work has gone into the computational problem of computing the
zeta function of a ``random'' curve over a finite field, in large part because
the question of determining the order of the Jacobian group (class group)
of such a curve stands in the way of using said group for public key
cryptography. The history of this problem is not our present concern, and
anyway it has been documented elsewhere; see for example
\cite{me-ants} for an overview.

By contrast, relatively little work has gone into the analogous computational
problem over higher dimensional varieties. Ongoing work of Bas Edixhoven 
and his
collaborators, to give an efficient algorithm for computing the $n$-th
Fourier coefficient of a fixed modular form when $n$ is a large integer
of known factorization, amounts to computing factors of zeta functions of
higher-dimensional varieties over large prime fields using methods of
$\ell$-adic cohomology. Over fields of small characteristic, one also may
use techniques of $p$-adic cohomology, which when applicable tend to yield
more efficient algorithms. However, while a number of 
reasonable-looking algorithms for various higher-dimensional varieties
have been described theoretically, e.g., by Gerkmann \cite{gerkmann},
Lauder \cite{lauder-def, lauder-def2}, 
and Lauder and Wan \cite{lauder-wan},
until recently nothing had been attempted in practice.
(For some very recent developments on this front, see
Section~\ref{subsec:alternate}.)

The purpose of this paper is to begin repairs on this gap in knowledge, by
on one hand illustrating how even limited information about the action of
Frobenius in $p$-adic cohomology can be used to address questions of some
possibly practical import, and on the other hand to outline an algorithm
which has been demonstrated in practice to be able to obtain this limited
information. The potential import stems from the fact that one can use
information about Frobenius, specifically bounds on Picard numbers obtained
from performing linear algebra on a low-precision Frobenius matrix,
to control the minimum distance of an algebraic
geometry (Goppa) code derived from a surface. As observed by Voloch and 
Zarzar,
such codes have the LDPC (low density parity check) property and so may be of
special interest.

Besides this introduction, the paper is structured in four main sections.
The first section is general, describing in detail what a Picard number is
and how to use an approximately
computed Frobenius matrix to bound it. The second part
gathers some facts about algebraic de Rham cohomology and $p$-adic cohomology.
The third part 
sketches a particular algorithm for producing an approximate Frobenius
matrix on the cohomology of a smooth hypersurface, using $p$-adic cohomology
and a description of the cohomology of a smooth hypersurface
due to Griffiths \cite{griffiths};
we also mention some related proposals.
The fourth part describes an implementation of our algorithm in
\textsc{Magma}
and tallies a few experimental results.

\section{Picard numbers and Frobenius matrices}

\subsection{Picard groups}

\begin{defn}
Let $X$ be a variety over a field $k$. The \emph{Picard group} 
$\Pic(X)$ is the group of isomorphism classes of line bundles (or invertible
sheaves) on $X$. Note that for $X$ smooth, isomorphism classes of 
line bundles are in natural bijection with rational equivalence classes of
(Weil or Cartier) divisors on $X$.
\end{defn}

\begin{lemma} \label{L:compare picard}
Let $X$ be a smooth proper irreducible
variety over a field $k$, let $k^{\sep}$ denote the
separable closure of $k$, and put $G = \Gal(k^{\sep}/k)$. Then the natural map
\[
\Pic(X) \to \Pic(X \times_k 
k^{\sep})^{G}
\]
(in which the superscripted $G$ means take $G$-invariants)
is always injective; moreover, if the Brauer group $\Br(k)$ is trivial
(e.g., if $k$ is finite), then the map is  surjective.
\end{lemma}
\begin{proof}
For any smooth $X$, the Hochschild-Serre spectral sequence in
\'etale cohomology gives rise to an exact sequence
\[
0 \to H^1(G, k^{\sep}[X]^*) \to \Pic(X) \to \Pic(X \times_k 
k^{\sep})^{G}
\to H^2(G, k^{\sep}[X]^*)
\]
where $k^{\sep}[X]^* = H^0(X \times_k k^{\sep}, \GG_m)$ (see, e.g.,
\cite[Lemme~6.3]{sansuc}).
For $X$ proper irreducible, $k^{\sep}[X]^* = (k^{\sep})^*$, and
$H^1(G, (k^{\sep})^*) = 0$ by Hilbert's Theorem 90, while
$H^2(G, (k^{\sep})^*) = \Br(k)$. 
This yields the desired results.
\end{proof}

\begin{defn}
For $X$ smooth proper irreducible 
over a field $k$, we say that a divisor is 
\emph{algebraically equivalent
to zero} if it has the form 
\[
(D \cap (X \times \{p\})) - (D \cap (X \times \{q\}))
\]
for some connected (but not necessarily smooth or irreducible) 
curve $C$, some pair of closed points $p,q$ on $C$ of the same degree, and some
divisor $D$ on $X \times C$ containing no fibres of the projection
$X \times C \to C$. The set of divisors algebraically equivalent to zero is a
subgroup closed under rational equivalence; let $\Pic^0(X)$ denote the image of
this subgroup in $\Pic(X)$.
\end{defn}

\begin{remark}
For $k$ algebraically closed,
the elements of $\Pic^0(X)$ can be thought of in a natural way as the
$k$-valued points on a certain variety over $k$, the \emph{Picard variety}
associated to $X$; there is also a scheme-theoretic version of this
fact that works over more general bases. 
We will not use this interpretation here.
\end{remark}

\begin{lemma} \label{L:torsion-free}
Let $X$ be a smooth irreducible complete intersection in $\PP^r_k$, for
$k$ a perfect field and $r$ a positive integer. Then
$\Pic^0(X) = 0$, 
$\Pic(X)$ is torsion-free, and $\calO(1)$ is indivisible in $\Pic(X)$.
\end{lemma}
\begin{proof}
For $k$ algebraically closed, this is 
\cite[Th\'eor\`eme~1.8]{deligne};
the general case follows from the algebraically closed case plus
Lemma~\ref{L:compare picard}. Note that if $\dim(X) \geq 3$, one in fact
has $\Pic(X) = \ZZ \cdot \calO(1)$.
\end{proof}

\subsection{N\'eron-Severi groups}

\begin{defn}
For $X$ smooth proper irreducible 
over a field $k$, the quotient $\NS(X) = \Pic(X)/\Pic^0(X)$ is called the
\emph{N\'eron-Severi group} of $X$. 
\end{defn}

\begin{remark}
For $X$ projective, one may define the degree of a divisor with respect to
any fixed ample divisor. The resulting map induces a 
homomorphism $\deg: \NS(X) \to \ZZ$
sending any ample divisor to a positive integer; in particular, $\NS(X)$ is 
nontrivial and any ample divisor represents a nonzero class in $\NS(X)$.
\end{remark}

\begin{remark}
 For $k=\CC$, there is a natural
map
\begin{equation} \label{eq:lefschetz11} 
\NS(X) \to H^2(X^{\an}, \ZZ) \cap H^{1,1}(X),
\end{equation}
where
$X^{\an}$ denotes the analytic space associated to $X$ and
$H^{1,1}(X) = H^1(X, \Omega^1_X)$ (which by GAGA may be computed either
algebraically or analytically); the Lefschetz (1,1)-theorem 
\cite[\S 1.2]{griffiths-harris}
asserts that the map \eqref{eq:lefschetz11} is a bijection.
\end{remark}

\begin{defn}
Let $X$ be a smooth proper irreducible variety over a field $k$.
From the Lefschetz theorem, it follows that $\NS(X)$ is finitely generated
in case $k = \CC$. In fact, $\NS(X)$ is always finitely generated;
this was first shown by N\'eron \cite[Th\'eor\`eme~2, p.~145]{neron} (see also 
\cite{hartshorne}, \cite[\S 1]{lang-neron}).
The rank of $\NS(X)$ as a finitely generated abelian group is called the
\emph{Picard number} (or \emph{arithmetic Picard number}) of $X$. 
The rank of $\NS(X \times_k \overline{k})$, for
$\overline{k}$ the algebraic closure of $k$, 
is called the \emph{geometric Picard
number} of $X$.
\end{defn}

\begin{defn}
Let $X$ be a smooth projective surface over a field $k$.
Then intersection theory \cite[Chapter~V]{hartshorne-book} gives rise to a
symmetric pairing on divisors of $X$, called the \emph{intersection pairing};
this pairing respects algebraic equivalence, so it descends to $\NS(X)$.
We say a divisor $D$ is \emph{numerically equivalent to zero} if 
$C \cdot D = 0$
for every projective curve $C$ contained in $X$; this turns out to happen
if and only if some multiple of $D$ is algebraically equivalent to zero
\cite[Theorem~4]{matsusaka}. That is,
the quotient of $\NS(X)$ by the classes of divisors numerically equivalent
to zero is precisely $\NS(X)/\NS(X)_{\tors}$. This group is never zero
because an ample divisor $D$ satisfies $C \cdot D > 0$
for any $C$, and so is not numerically equivalent to zero. That is, the
Picard number of a smooth projective surface is always positive.
\end{defn}

\begin{remark}
Although we will only attempt to bound Picard numbers over finite fields,
doing so also has consequences over number fields. That is because if
$\gothp$ is a prime ideal in the ring of integers $\gotho_K$
of a number field $K$,
$k = \gotho_K/\gothp$, and 
$X$ is a smooth projective surface over the localization
of $\gotho_K$ at $\gothp$, then the torsion-free quotient of $\NS(X_K)$
injects into the torsion-free quotient of $\NS(X_k)$, compatibly with the
intersection pairing \cite[\S 6]{vanluijk-heron}. We can thus control the
size of $\NS(X_K)$ by controlling $\NS(X_k)$; in some cases, one can gain
further control by reducing modulo more than one prime of good reduction
\cite{vanluijk-pic}.
\end{remark}

\subsection{Picard numbers and codes}

We now recall briefly what Picard numbers 
have to do with error-correcting codes;
the link lies in a
higher-dimensional version of Goppa's construction
\cite{goppa} of algebraic geometry codes from curves over finite fields.

\begin{defn}
Let $X$ be a smooth 
projective irreducible variety over a finite field $\FF_q$. Let $H$ be
an ample divisor on $X$, let $m$ be a positive integer such that the divisor
$mH$ is very ample, and let $\calL(mH) = \Gamma(X, \calO(mH))$ be the Riemann-Roch
space of $mH$; we may identify elements of $\calL(mH)$ with rational
functions $f \in \FF_q(X)$ such that the divisor $\divi(f) + mH$ is effective.
Let $S$ be the set of $\FF_q$-rational points of $X \setminus H$.
Define the code $C(X,mH)$ 
to be the subspace of $\FF_q^S$ of functions induced by elements
of $\calL(mH)$, viewed as a linear code over $\FF_q$.
\end{defn}

In Goppa's original construction, $X$ is a curve, the rate of the code
(the ratio between the dimensions of the code and of its ambient space)
is determined by Riemann-Roch, and a good bound on the minimum distance 
(the smallest number of nonzero elements in a nonzero codeword) comes
from the fact that a rational function cannot have more zeros than poles.
In higher dimensions, one can still get rate information out of 
Riemann-Roch, but bounding the minimum distance is trickier.
For surfaces, this question has been investigated by Voloch and Zarzar, with
the idea of using the subcodes induced by curves on a surface to give
an asynchronous decoding algorithm in the style of Luby-Mitzenmacher
\cite{luby-mitzenmacher}.
Voloch and Zarzar observe that a low Picard number gives rise to a good
bound on the minimum distance; we limit ourselves here to mentioning two
sample assertions, and defer to 
\cite{voloch-zarzar} for more information.

\begin{lemma}[Voloch]
Let $X$ be a smooth projective surface over a field $k$, and suppose 
$\NS(X)/\NS_{\tors}(X)$ 
is generated by the ample divisor $H$. Then for any positive integer
$m$, the zero divisor of a nonzero element of $\calL(mH)$ has at most
$m$ irreducible components.
\end{lemma}
\begin{proof}
See \cite[Lemma~2.2]{voloch-zarzar}.
\end{proof}

\begin{lemma}[Zarzar]
Let $X$ be a smooth surface of degree $d$ in $\PP^3$ over a perfect field $k$,
and suppose that the Picard number of $X$ is equal to $1$. If $Y$ is
an irreducible surface in $\PP^3$ of degree $m<d$, then
$X \cap Y$ is also irreducible.
\end{lemma}
\begin{proof}
First note that $\Pic(X)$ and $\NS(X)$ coincide and are torsion-free by
Lemma~\ref{L:torsion-free}. Then invoke
\cite[Lemma~2.1]{zarzar}.
\end{proof}

\subsection{Zeta functions and Picard numbers}
\label{sec:zeta}

\begin{defn} \label{D:zeta}
Let $X$ be a smooth proper
variety over a finite field $\FF_q$. The \emph{zeta function}
of $X$ is the power series
\[
Z(X,T) = \exp \left( \sum_{n=1}^\infty \#X(\FF_{q^n}) \frac{T^n}{n} \right).
\]
The Weil conjectures, proved by Dwork, Grothendieck, Deligne, et al.\
(see \cite[Appendix~C]{hartshorne-book} for a fuller statement), assert
that there exists a product decomposition
\[
Z(X,T) = \prod_{i=0}^{2 \dim(X)} P_i(T)^{(-1)^{i+1}},
\]
where $P_i(T) \in \ZZ[T]$ and $P_i(0) = 1$, such that the roots of 
$P_i(T)$ in $\CC$ all have absolute value $q^{-i/2}$. Moreover, the
multiset of roots of $P_i$ is invariant under the transformation
$r \mapsto q^{-i}/r$.
\end{defn}

The connection between zeta functions and Picard numbers was first
articulated by Tate \cite{tate}, who showed that
\begin{equation} \label{eq:tate}
\rank \NS(X) \leq \ord_{T = 1/q} P_2(T).
\end{equation}
(Actually Tate's argument gives a bit more information than this; see
Remark~\ref{rem:tate2} below.)
Tate conjectured further (by analogy with the conjecture of Birch
and Swinnerton-Dyer) that equality always holds in \eqref{eq:tate};
Tate himself showed this for abelian varieties, 
and it is also known in some other cases, e.g., for ordinary K3 surfaces
\cite{zarhin}. Tate's conjecture in general is wide open; however, since our
purpose here is merely to give upper bounds for Picard numbers, the 
unconditional
inequality \eqref{eq:tate} will suffice.
(See \cite{tate2} for more context on Tate's conjecture.)

Note that \eqref{eq:tate} also gives a bound on the geometric Picard number:
\begin{equation} \label{eq:tate2}
\rank \NS(X \times_{\FF_q} \overline{\FF_q}) \leq \sum_\zeta 
\ord_{T = \zeta/q} P_2(T),
\end{equation}
where $\zeta$ runs over all roots of unity. 
Since $P_2(T)$ has integer coefficients,
we can rewrite \eqref{eq:tate2} as
\[
\rank \NS(X \times_{\FF_q} \overline{\FF_q}) \leq \sum_{n: \phi(n) \leq \deg(P_2)} 
\phi(n) \ord_{T = \zeta_n/q}
P_2(T),
\]
where $\zeta_n$ denotes any one primitive $n$-th root of unity.
For computational purposes, we need an explicit bound on $n$;
here is an easy such bound.
\begin{lemma} \label{L:phibound}
For any positive integer $n$, we have 
\[
\phi(n) \geq \frac{n}{\lfloor \log_2(n) \rfloor + 1}.
\]
\end{lemma}
\begin{proof}
We have 
\[
\frac{\phi(n)}{n} = \prod_p \left( 1 - \frac{1}{p} \right),
\]
the product running over the distinct prime divisors of $n$. There are at most
$\lfloor \log_2(n) \rfloor$ such divisors, so
\[
\frac{\phi(n)}{n} \geq \prod_{i=2}^{\lfloor \log_2(n) \rfloor + 1}
\left( 1 - \frac{1}{i} \right) = \frac{1}{\lfloor \log_2(n) \rfloor + 1},
\]
as desired.
\end{proof}

Here is a standard parity consideration that arises in the context
of Tate's conjecture.

\begin{remark} \label{rem:rank2}
Note that the right side of \eqref{eq:tate2} has the same parity as 
$\ord_{T=1/q} P_2(T) + \ord_{T=-1/q} P_2(T)$. This in turn has the same
parity as $\deg(P_2)$, since $\pm 1/q$ are the only real roots consistent
with the restriction that each root has absolute value $1/q$.
Under Tate's conjecture, equality would hold in \eqref{eq:tate2}, and would
thus imply that if $X$ is a smooth proper variety
 for which $\deg(P_2)$ is even, then the 
geometric Picard number of $X$ is at least 2.
\end{remark}

\subsection{Weil cohomologies}

\begin{defn}
Fix a finite field $\FF_q$ of characteristic $p$, 
and let $K$ be a field of characteristic zero.
By a \emph{(weak) Weil cohomology} over $K$, we will mean the following data.
\begin{itemize}
\item
One must specify a collection of
contravariant functors from smooth proper 
varieties $X$ over $\FF_q$ to finite dimensional $K$-vector spaces 
$H^i(X)$ equipped with endomorphisms $F_i$, such that
\[
P_i(T) = \det(1 - TF_i, H^i(X)) \qquad (i=0, \dots, 2 \dim(X)).
\]
For $m \in \ZZ$, we write $H^i(X)(m)$ to mean the vector space $H^i(X)$ equipped with
the endomorphism $q^{-m} F_i$.
\item
The Lefschetz trace formula holds for Frobenius: for any positive
integer $n$,
\[
\#X(\FF_{q^n}) = \sum_{i=0}^{2 \dim(X)} (-1)^i \trace(F_i^n, H^i(X)).
\]
\item
One must specify functorial, $F$-equivariant maps (for $d = \dim(X)$)
\[
\trace_X: H^{2d}(X)(d) \to K
\]
(where $F$ acts as the identity on $K$)
which should be isomorphisms when $X$ is geometrically irreducible.
\item
One must specify associative, functorial, $F$-equivariant
cup product pairings $\cup: H^i(X) \times H^j(X) \to H^{i+j}(X)$ 
such that (for $d = \dim(X)$)
the pairings
\[
H^i(X) \times H^{2d-i}(X)(d)
\stackrel{\cup}{\to} H^{2d}(X)(d) 
\stackrel{\trace_X}{\to} K
\]
are perfect (Poincar\'e duality).
\item
One must specify an injective $K$-linear homomorphism
\[
\NS(X) \otimes_{\ZZ} K \to H^2(X)(1)^{F=1}
\]
(the cycle class map).
\end{itemize}
For a more precise definition of the phrase ``Weil cohomology''
(which actually includes more structure than this, including
a K\"unneth decomposition, cycle class maps for higher Chow groups,
and a full Lefschetz hyperplane theorem, plus additional compatibilities),
see \cite{kleiman}.
\end{defn}

\begin{remark} \label{rem:tate2}
By virtue of the cycle class map, 
the existence of a Weil cohomology 
yields the inequality \eqref{eq:tate}, as the right side of 
\eqref{eq:tate} equals the dimension of the generalized eigenspace of
$H^2(X)$ with eigenvalue $q$. 
In practice, we will use the resulting
slightly stronger version of \eqref{eq:tate}:
\begin{equation} \label{eq:tate3}
\rank \NS(X) \leq \corank(F_2-q, H^2(X))
\end{equation}
and the corresponding version of \eqref{eq:tate2}:
\begin{equation} \label{eq:tate4}
\rank \NS(X \times_{\FF_q} \overline{\FF_q}) \leq \sum_{n: \phi(n) \leq \deg(P_2)} 
\phi(n) \corank(F_2 - \zeta_n q, H^2(X)).
\end{equation}
In theory, one expects that $F_2$ acts semisimply on $H^2(X)$; this would
follow from the full conjecture of Tate, which is somewhat stronger
than we have described here (as it makes predictions about $H^{2i}(X)$ for
all $i$).
\end{remark}

\begin{remark}
At the time \cite{tate} was written, the only Weil cohomologies that had
been constructed were the $\ell$-adic \'etale cohomologies for each
$\ell \neq p$, which take values in $\QQ_{\ell}$. 
Subsequently, it was shown that Berthelot's rigid cohomology is also a
Weil cohomology; it takes values in 
the $p$-adic field $\QQ_q$.
(Here and throughout, for brevity we write $\ZZ_q$ for $W(\FF_q)$ and
$\QQ_q$ for $\Frac \ZZ_q$.)
For the essential properties of rigid cohomology, see
\cite{berthelot1, berthelot2}; also see \cite{illusie} for additional
context on $p$-adic cohomology theories.
\end{remark}

\subsection{Approximate Gaussian elimination}

Let us set some notation for this subsection.
\begin{notation}
Let $K$ be a complete discretely valued field.
Let $\gotho_K$ be the
ring of integers of $K$, fix a uniformizer
$\pi$ of $\gotho_K$, and write $v(x)$ for the valuation of $x \in K$.
\end{notation}

We are going to describe an algorithm for producing an upper bound on the
corank of a matrix $A$ over $K$, 
given only the information of the entries of $A$
modulo $\pi^m$ for some integer $m$. There is no harm in rescaling the matrix
$A$
(by multiplying by an appropriate power $\pi^n$ of $\pi$, then replacing $m$
by $m+n$) to ensure that $A$ has entries in $\gotho_K$.

\begin{algo} \label{algo:picard}
Given a matrix $\overline{A}$ with entries in $\gotho_K / \pi^m \gotho_K$
which is the reduction of a matrix $A$ over $\gotho_K$, the following algorithm
returns an upper bound on $\corank(A)$.
\begin{enumerate}
\item
Let $r$ be the number of rows of $\overline{A}$, and let
$s$ be the number of  columns of $\overline{A}$.
If $\overline{A}$ has no nonzero entries (possibly because $\overline{A}$
is an empty matrix), return $s$ and STOP.
\item
Choose a nonzero entry $\overline{A}_{ij}$ of minimum valuation.
\item
For $k = 1, \dots, r$ in succession, skipping over $k=i$,
choose $c \in \gotho_K/\pi^n \gotho_K$
such that $c \overline{A}_{ij} = \overline{A}_{kj}$, and subtract 
$c$ times the $i$-th row of $\overline{A}$ from the $k$-th row of
$\overline{A}$.
\item
For $\ell = 1, \dots, s$ in succession, skipping over $\ell=j$,
choose $c \in \gotho_K/\pi^n \gotho_K$
such that $c \overline{A}_{ij} = \overline{A}_{i\ell}$, and subtract 
$c$ times the $j$-th column of $\overline{A}$ from the $\ell$-th column of
$\overline{A}$.
\item
Delete row $i$ and column $j$ from $\overline{A}$, then go to step 1.
\end{enumerate}
\end{algo}
\begin{proof}
We prove the claim by induction on the number of rows plus columns of $A$.
If $\overline{A}$ is the zero matrix, the claim is evident. Otherwise,
we may assume without loss of generality that $i=j=1$. By lifting each
of the row and column operations from $\overline{A}$ to $A$ appropriately, 
we may also assume that $A_{i1} = 0$ for $i=2,\dots, m$ and that
$A_{1j} = 0$ for $j=2,\dots, n$. It is now clear that the corank of $A$
is equal to that of its lower right $(m-1)\times(n-1)$ submatrix.
\end{proof}

One can also use Algorithm~\ref{algo:picard} 
to obtain information about the determinant
of $A$.
\begin{prop} \label{prop:det}
Given an $n \times n$
matrix $\overline{A}$ with entries in $\gotho_K / \pi^m \gotho_K$
which is the reduction of a matrix $A$ over $\gotho_K$, suppose that
Algorithm~\ref{algo:picard} returns the bound $0$. For $h = 1, \dots,n$,
with  notation as in the $h$-th iteration of the algorithm, put
$\overline{a_h} = \overline{A}_{ij} \in \gotho_K/\pi^m \gotho_K$ and $e_h = (-1)^{i+j}
\in \gotho_K$.
Choose lifts $a_1, \dots, a_n$ of 
$\overline{a_1}, \dots, \overline{a_n}$ to $A$. Then
\begin{equation} \label{eq:det bound}
v(\det(A) - a_1e_1\cdots a_ne_n) \geq
\min_i \{m - v(a_i)\} + \sum_{i=1}^n v(a_i).
\end{equation}
\end{prop}
\begin{proof}
Perform the ``shadow'' computation of
the proof of Algorithm~\ref{algo:picard} with the following change: at each
step, instead of deleting row $i$ and column $j$, move them to the far
bottom and right. The final matrix has determinant 
$e_1\cdots e_n\det(A)$; on the other hand, it is diagonal with
entries congruent to $a_1, \dots, a_n$ modulo $\pi^m$.

To obtain the desired estimate, expand $\det(A)$ as a sum of
signed products $\pm A_{1\sigma(1)} \cdots A_{n \sigma(n)}$
with $\sigma$ a permutation.
When measuring the valuation of the product,
each term $A_{ij}$ with $i=j$ contributes $v(a_i)$ and each
term $A_{ij}$ with $i \neq j$ contributes $m$; hence for $\sigma$
different from the identity, this product
has valuation at least the right side of \eqref{eq:det bound}.
The product $A_{11} \cdots A_{nn} - a_1 e_1 \cdots a_n e_n$
is the sum of the expressions
\[
a_1 e_1 \cdots a_{j-1} e_{j-1} A_{jj} \cdots A_{nn}
- a_1 e_1 \cdots a_j e_j A_{(j+1)(j+1)} \cdots A_{nn}
\]
for $j=1, \dots, n$, which have respective valuations
$m - v(a_j) + \sum_{i=1}^n v(a_i)$. This yields the claim.
\end{proof}

\begin{remark} \label{R:char poly}
Note that Proposition~\ref{prop:det} may be used to obtain approximate
characteristic polynomials, by applying it with the field $K$ replaced by
the completion of $K(t)$ for the Gauss valuation (i.e., the valuation which
on a polynomial returns the minimum valuation of any coefficient)
and approximating $\det(tI - A)$.
If one knows that the matrix $A$ is ``nearly'' divisible by some $\pi^i$,
then one may obtain better information by approximating 
$\det(\pi^i t I - A)$; this often happens in the setting of
$p$-adic cohomology.
For example, suppose $X$ is a smooth proper variety over $\ZZ_q$,
where $q$ is a power of a prime $p < \dim(X_{\FF_q})$. Then if one writes
down the Frobenius matrix on the $j$-th rigid 
cohomology with respect to an appropriate
basis (namely, a basis of crystalline cohomology modulo torsion), 
the Hodge numbers 
\[
h^{i,j-i} = \dim_{\QQ_q} H^{j-i}(X_{\QQ_q}, \Omega^i_{X/\QQ_q}) \qquad
(i=0,\dots,j)
\]
give the multiplicities of $p^i$ as elementary divisors of the matrix.
See \cite[Theorem~1.3.9]{illusie} for a more general statement.
\end{remark}

\begin{remark}
It may be possible to obtain even better bounds on characteristic polynomials
which are more adaptive (i.e., give individual bounds for each coefficient)
by using more careful linear algebra plus $p$-adic floating point arithmetic.
(In fact, the bound in Proposition~\ref{prop:det} is most naturally
phrased in terms of floating point arithmetic: the first term in
\eqref{eq:det bound} is the minimum precision of the mantissa of an entry
in the final matrix.)
In particular, it would be interesting
to give such bounds
for the more general setting where the accuracy may vary from
entry to entry; in our application to bounding Picard numbers,
being able to work in this generality 
might lend some flexibility to the cohomological calculation. We will
not consider the more general setting here.
\end{remark}

\section{A little $p$-adic cohomology}

In this part, we set up a bit of the theory of $p$-adic cohomology
for use later on; this involves some consideration of
 algebraic de Rham cohomology. 
Some of the calculations, particularly 
Theorem~\ref{thm:deligne}, may be of independent interest.

\subsection{A homological calculation}

We start with a brief excursion into homological algebra, following
\cite[Chapter~17]{eisenbud}.

\begin{defn} \label{D:koszul}
Let $R$ be a ring and choose $x = (x_1, \dots, x_n) \in R^n$.
The \emph{Koszul complex} $K(x)$ is the exterior algebra
$\wedge^*_R (R^n)$ with differentials
\[
d_x(z) = x \wedge z.
\]
Let $K'(x)$ denote the dual complex, whose underlying $R$-module
we also identify with $\wedge^*_R (R^n)$; let
$\del_x$ denote the differentials in $K'(x)$.
\end{defn}

\begin{lemma} \label{L:homotopy}
Let $R$ be a ring. For any $x,y \in R^n$ and any $z \in \wedge^*_R (R^n)$,
\[
(d_x \del_y + \del_y d_x)(z) = (x_1y_1 + \cdots + x_ny_n)z.
\]
\end{lemma}
\begin{proof}
An easy calculation: see \cite[Lemma~17.13]{eisenbud}.
\end{proof}

\begin{prop} \label{P:homotopy}
Let $R$ be a ring, and let
$C$ be a complex in the category of $R$-modules.
Then for $(x_1, \dots, x_n) \in R^n$,
the homology of the product
complex $K(x) \otimes C$ is annihilated by the ideal
$(x_1, \dots, x_n)$.
\end{prop}
\begin{proof}
Given $r$ in the ideal $(x_1, \dots, x_n)$, choose $y_1, \dots, y_n \in R$ such that
$x_1 y_1 + \cdots + x_n y_n = r$. Then Lemma~\ref{L:homotopy} shows that
multiplication by $r$ is homotopic to zero on $K(x)$, with the
homotopy given by $\del_y$. Tensoring that homotopy with the identity map
on $C$ yields the same assertion on $K(x) \otimes C$, proving the claim.
\end{proof}

\subsection{Algebraic de Rham cohomology}

Algebraic de Rham cohomology is usually considered over a field of
characteristic zero, but for $p$-adic cohomological calculations, we
also need to work with it over arithmetically interesting base schemes.

\begin{defn}
By a \emph{smooth pair} (resp.\ \emph{smooth proper pair})
of relative dimension $n$ over a scheme $S$,
we mean a pair $(X,Z)$ in which $X$ is a smooth (resp.\ smooth proper)
$S$-scheme of relative dimension $n$, 
and $Z$ is a relative reduced normal crossings
divisor on $X$. That is (in the smooth-only case), 
\'etale locally on $X$, we can find an 
$S$-isomorphism
of $X$ with a Zariski open subset of $\AAA^n_S$ under which $Z$ is carried
to an open subset of a union of some 
(or all, or none) of the coordinate hyperplanes.
We think of $Z$ as defining a logarithmic structure on
$X$ in the sense of Kato \cite{kato}, and write $(X,Z)$ also for the resulting
log scheme. If $Z = \emptyset$, we abbreviate $(X,Z)$ to simply $X$.
\end{defn}

\begin{defn}
Let $(X,Z)$ be a smooth pair of relative dimension $n$ over a scheme $S$.
Let $\Omega_{(X,Z)/S}$ be the sheaf of differentials on $X$ with logarithmic
singularities along $Z$ (i.e., one allows $dt/t$
for $t$ a local parameter of a component of $Z$), 
relative to $S$; then $\Omega_{(X,Z)/S}$ is a locally
free coherent $\calO_X$-module of rank $n$. 
Put $\Omega^i_{(X,Z)/S} = \wedge^i_{\calO_X} \Omega_{(X,Z)/S}$; then exterior
differentiation induces maps $d_i: \Omega^i_{(X,Z)/S} \to 
\Omega^{i+1}_{(X,Z)/S}$
such that $d_{i+1} \circ d_i = 0$. The resulting complex is called the
\emph{de Rham complex} of $(X,Z)$ over $S$, and its $j$-th hypercohomology 
\[
H^j_{\dR}((X,Z)/S) = \HH^j(X, \Omega^._{(X,Z)/S})
\]
is called the \emph{$j$-th algebraic de Rham cohomology group}
of $(X,Z)$ over $S$.
\end{defn}

\begin{remark}
Since the $\Omega^i$ are quasi-coherent, we may calculate algebraic de Rham
cohomology on the \'etale site instead of the Zariski site
\cite[Proposition~3.7]{milne}.
 This permits the use of \'etale localization arguments.
\end{remark}


Let $(X,Z)$ be a smooth pair over $\CC$,
and put $U = X \setminus Z$; one then has an isomorphism
\begin{equation} \label{eq:log sing}
H^j_{\dR}((X,Z)/\CC) \cong H^j_{\dR}(U).
\end{equation}
Namely, by Serre's GAGA theorem 
\cite{serre} on the left side and Grothendieck's comparison theorem
\cite{grothendieck} on the right side (which also uses GAGA, together
with resolution of singularities), this may be checked for 
analytic de Rham cohomology, where it amounts to the
Poincar\'e lemma (see \cite{deligne} for variations on this theme).

However, it was pointed out to us by Johan de Jong that one can also
establish \eqref{eq:log sing} algebraically.
In so doing, one can also prove an integral variant where one compares
cohomology of the complex of differentials with logarithmic poles
with the complex of differentials where the poles are made somewhat worse.
Here is the result; in its relevance to computing
$p$-adic cohomology, it should
be viewed as a generalization of 
\cite[Lemma~2]{me-count}.

\begin{defn}
Let $(X,Z)$ be a smooth pair over a scheme $S$, and put 
$U = X \setminus Z$. For $m$ a nonnegative integer, write (as usual)
$\Omega^j_{(X,Z)/S}(mZ)$ for the twist $\Omega^j_{(X,Z)/S} \otimes_{\calO_X} \calO_X(mZ)$;
note that $\Omega^j_{U/S}$ is the direct limit of the
$\Omega^j_{(X,Z)/S}(mZ)$ as $m$ increases.
\end{defn}

\begin{theorem} \label{thm:deligne}
Let $(X,Z)$ be a smooth pair over a scheme $S$.
For each nonnegative integer $m$,
the cokernels of the maps on homology sheaves
induced by the natural map of complexes of sheaves
\[
\Omega^._{(X,Z)/S} \to \Omega^._{(X,Z)/S}(mZ)
\]
are killed by $\lcm(1, \dots, m)$.
\end{theorem}
\begin{proof}
The claim may be verified stalkwise on $S$, so we may assume $S = \Spec A$
is affine and local.
We may also localize on $X$; starting with a point on $X$, we can shrink
to ensure that $X = \Spec R$ is affine,
$\Omega^._{X/S}$ is generated
freely by $dx_1, \dots, dx_n$ for some $x_1, \dots, x_n \in R$,
$Z = V(x_1 \cdots x_h)$, and our chosen point lies on $V(x_1, \dots,x_h)$.
For each $T \subseteq \{1, \dots, h\}$, let $I_T$ be the ideal
generated by $x_i$ for each $i \in T$, and put $R_T = R/I_T$.
By a further Zariski localization, we may also assume that
$R$ contains a copy of each $R_T$.

Define
\[
\tilde{d}x_i = \begin{cases} \frac{dx_i}{x_i} & 1 \leq i \leq h \\
dx_i & i > h
	       \end{cases}
\qquad
\tilde{\del}_i = \begin{cases} x_i \frac{\del}{\del x_i} & 1 \leq i \leq h \\
\frac{\del}{\del x_i} & i > h.
		 \end{cases}
\]
For each subset $U = \{i_1, \dots, i_r\}$ of $\{1, \dots, n\}$ with
$i_1 < \cdots < i_r$, put 
\[
\tilde{d} x_U = \tilde{d}x_{i_1} \wedge \cdots \wedge \tilde{d}x_{i_r}.
\]

Let $M$ be the set of monomials $\mu = x_1^{j_1} \cdots x_h^{j_h}$ with
$j_1, \dots, j_h \in \{0, \dots, m\}$, viewed as a partially ordered 
set under divisibility. For $D$ a nonempty down-closed subset of $M$
(i.e., one in which inclusion of $\mu$ implies inclusion of any divisor
of $\mu$, so that in particular $1 \in D$), define
\[
Q_D = \bigcup_{\mu \in D} \mu^{-1} R;
\]
note that each $\tilde{\del}_i$ sends $Q_D$ into itself,
so that $Q_D \otimes_R \Omega^._{(X,Z)/S}$ is again a complex.

Let $D$ be a down-closed subset of $M$ strictly bigger than $\{1\}$.
Choose $\mu = x_1^{j_1} \cdots x_h^{j_h}$ maximal in $D$; then
$D' = D \setminus \{\mu\}$ is also down-closed.
Let $T$ be the set of $i \in \{1, \dots, h\}$ for which $j_i \neq 0$,
which is necessarily nonempty since we cannot have $\mu = 1$;
then the chosen inclusion $R_T \hookrightarrow R$ induces an isomorphism
\[
\mu^{-1} R_T \cong Q_D/Q_{D'}.
\]
This isomorphism is equivariant for the action of the $\tilde{\del}_i$ for
$i \notin T$; for $i \in T$, it converts multiplication by $-j_i$ on the left
side into the induced action of
$\tilde{\del}_i$ on the right. It consequently induces an 
isomorphism of complexes
\begin{align*}
\ZZ(j_i: i \in T) \otimes_{\ZZ} \Omega^._{R_T/A}
&\cong (Q_D/Q_{D'}) \otimes_R \Omega^._{(X,Z)/S} \\
&\cong (Q_D \otimes_R \Omega^._{(X,Z)/S})/(Q_{D'} \otimes_R \Omega^._{(X,Z)/S}).
\end{align*}
By Proposition~\ref{P:homotopy}, the homology of these complexes
is killed by $\gcd(j_i: i \in T)$.

Now 
suppose $\omega \in Q_D \otimes_R \Omega^r_{(X,Z)/S}$ satisfies $d\omega = 0$.
Write $\omega = \sum_U g_U \tilde{d}x_U$,
where $U$ runs over $r$-element subsets of $\{1, \dots, n\}$.
Let $h_U \in \mu^{-1} R_T$ be the image of $g_U$ under the map
$Q_D/Q_{D'} \to \mu^{-1} R_T$; then $d (\sum_U h_U \tilde{d}x_U) = 0$
as well. (It may clarify to think of $h_U$ as a ``leading term''
in an expansion of $g_U$ as a formal Laurent series in the $x_i$ for $i \in T$.)
By the previous paragraph, 
$\sum_U h_U \tilde{d}x_U$ times 
$\gcd(j_i: i \in T)$ is a coboundary.
Since $1 \leq \gcd(j_1: i \in T) \leq m$,
$\omega$ is thus equal to a differential
whose image in the cokernel of the map on 
homology is killed by $\lcm(1, \dots, m)$ plus
a cocycle in $Q_{D'} \otimes_R \Omega^r_{(X,Z)/S}$.
This yields the claim by induction on the cardinality of $D$.
\end{proof}
\begin{cor} \label{C:deligne}
Let $(X,Z)$ be a smooth pair over a field $K$ of characteristic $0$, and put
$U = X \setminus Z$. Then
for each $i$, the map
$H^i_{\dR}((X,Z)/K) \to H^i_{\dR}(U/K)$ is an isomorphism of $K$-vector spaces.
\end{cor}
\begin{proof}
Apply Theorem~\ref{thm:deligne} and
note that formation of cohomology commutes with direct limits.
\end{proof}

\begin{remark} \label{R:lower pole order}
Note that if $Z$ is smooth, then the proof of Theorem~\ref{thm:deligne}
actually gives something slightly stronger: any closed $r$-form with poles
of order $m+1$ can locally be written as a closed logarithmic $r$-form 
plus a form which when multiplied by $\lcm(1, \dots, m)$ becomes
the differential of an $(r-1)$-form with poles of order $\leq m$.
This is crucial for our application: see Proposition~\ref{P:prec bound1}.
\end{remark}

\begin{prop} \label{P:excision}
Let $(X,Z)$ be a smooth pair over a scheme $S$, with $Z$ also smooth;
let $j$ denote the inclusion $Z \hookrightarrow X$.
Then there is an exact sequence of complexes of coherent sheaves on $X$:
\[
0 \to \Omega^\cdot_{X/S} \to \Omega^\cdot_{(X,Z)/S} \stackrel{\Res}{\to}
j_* \Omega^{\cdot}_{Z/S}[+1] \to 0,
\]
yielding an exact sequence in cohomology
\[
\cdots \to       H^i_{\dR}(X/S) \to H^i_{\dR}((X,Z)/S) \to H^{i-1}_{\dR}(Z/S)
\to H^{i+1}_{\dR}(X/S) \to \cdots.
\]
\end{prop}
\begin{proof}
Locally on $X$ we can choose a coordinate $x_n$ cutting out $Z$; then
each section of the quotient 
$\Omega^r_{(X,Z)/S}/\Omega^r_{X/S}$
admits a representative of the form
$dx_n/x_n \wedge \omega$ for some $\omega$, and the residue
map carries this section to the reduction of $\omega$ modulo $x_n$. To see that
this is well-defined globally, we must simply observe that the map is not
changed by changing the choice of the parameter $x_n$: if $u$ is a unit on
$X$, then $d(x_n u)/(x_n u) = dx_n/x_n + du/u$,
so $dx_n/x_n \wedge \omega$ and $d(x_n u)/(x_n u) \wedge \omega$ represent the
same element of the quotient. 
\end{proof}

\subsection{$p$-adic cohomology in general}

We suggest \cite{illusie} as a useful survey for the material underlying this
subsection.

\begin{notation}
Throughout this subsection, let $k$ be a perfect field of characteristic $p$,
and write $K$ for $\Frac W(k)$ and $\gotho_K$ for $W(k)$.
Let $\sigma: \gotho_K \to \gotho_K$ denote the Witt vector Frobenius, which is
the unique lift to $\gotho_K$ of the absolute Frobenius on $k$.
\end{notation}

\begin{defn}
Let $(X,Z)$ be a smooth proper pair over $k$.
Let $H^i_{\crys}(X,Z)$ denote the $i$-th 
\emph{(log\nobreakdash-)crystalline cohomology} of $(X,Z)$; this is an
 $\gotho_K$-module whose
construction is contravariantly functorial in the pair $(X,Z)$.
Moreover, the absolute Frobenius on $k$ acts $\sigma$-semilinearly on 
$H^i_{\crys}(X,Z)$.
For the general construction, see for instance \cite[Chapter~2]{shiho}.
\end{defn}

To make the Frobenius action on de Rham cohomology explicit, we need to pass
to rigid cohomology, so we can use the Monsky-Washnitzer interpretation of
$p$-adic cohomology.

\begin{defn} \label{D:rigid}
Let $X$ be a smooth $k$-scheme. Let $H^i_{\rig}(X)$ denote the 
\emph{$i$-th rigid cohomology group} of $X$ in the sense of Berthelot
\cite{berthelot1}; it is a finite-dimensional $K$-vector space whose construction
is contravariantly functorial in $X$.
Moreover, the absolute Frobenius on $k$ acts $\sigma$-semilinearly on 
$H^i_{\rig}(X)$.
For $X$ proper, Berthelot \cite[Proposition~1.9]{berthelot1} constructs
a functorial, Frobenius-equivariant isomorphism
\[
H^i_{\crys}(X) \otimes_{\gotho_K} K \cong H^i_{\rig}(X).
\]
For $(X,Z)$ a smooth pair over $k$, Shiho \cite[\S 2.4]{shiho} (using crucially
a result of Baldassarri and Chiarellotto \cite{bc}) constructs a functorial, 
Frobenius-equivariant
isomorphism
\[
H^i_{\crys}(X,Z) \otimes_{\gotho_K} K \cong H^i_{\rig}(X \setminus Z).
\]
\end{defn}

\subsection{$p$-adic cohomology in the liftable case}

When things can be lifted nicely to characteristic zero, the construction
of $p$-adic cohomology becomes much simpler.
\begin{prop} \label{P:liftable}
Let $(\frakX,\frakZ)$ be a smooth proper pair over $\gotho_K$.
Put $X = \frakX_k$, $Z = \frakZ_k$, and $U = X \setminus Z$.
\begin{enumerate}
\item[(a)]
There are isomorphisms
\[
H^i_{\crys}(X,Z) \cong H^i_{\dR}((\frakX,\frakZ)/\gotho_K)
\]
which are functorial in smooth proper pairs over $\gotho_K$. In particular, the right
side inherits an action of Frobenius.
\item[(b)]
Suppose that $\frakZ$ is also smooth, and use Proposition~\ref{P:excision},
Definition~\ref{D:rigid}, and (a) to obtain an exact sequence
\[
\cdots \to H^{i}_{\rig}(X) \to H^{i}_{\rig}(U) \to H^{i-1}_{\rig}(Z) \to 
H^{i+1}_{\rig}(X) \to \cdots.
\]
These maps are then Frobenius-equivariant for the following twists:
\begin{align*}
H^i_{\rig}(X) &\to H^i_{\rig}(U) \\
H^i_{\rig}(U) &\to H^{i-1}_{\rig}(Z)(-1) \\
H^{i-1}_{\rig}(Z) &\to H^{i+1}_{\rig}(X)(1),
\end{align*}
where $M(n)$ denotes $M$ with its absolute
Frobenius action multiplied by $p^{-n}$.
\end{enumerate}
\end{prop}
\begin{proof}
For (a), see \cite[Theorem~6.4]{kato}. For (b), we may invoke rigid analytic
GAGA to argue that algebraic de Rham cohomology of a 
smooth proper pair over $K$
coincides with rigid analytic de Rham cohomology of the 
analytification. The Frobenius equivariance may now be checked at the level of
complexes by following the construction of the Gysin isomorphism
in rigid cohomology \cite[\S 5]{berthelot1}, \cite{tsuzuki}.
\end{proof}

One also has a nice description of the Frobenius action on the rigid cohomology
of a smooth affine scheme in terms of a lifting.
\begin{defn} \label{D:monsky}
Let $X$ be a smooth affine $k$-scheme,
suppose $\frakX$ is a smooth affine $\gotho_K$-scheme lifting $X$,
write $\frakX = \Spec A$, 
and  choose a presentation $A \cong \gotho_K[x_1, \dots, x_n]/I$.
Let $\gotho_K \langle x_1, \dots, x_n \rangle^\dagger$ be the ring of
power series in $\gotho_K \llbracket x_1, \dots, x_n \rrbracket$ which
converge on some polydisc of radius greater than 1, and put
\[
A^\dagger = \gotho_K \langle x_1, \dots, x_n \rangle^\dagger/
I \gotho_K \langle x_1, \dots, x_n \rangle^\dagger.
\]
Then Berthelot \cite[Proposition~1.10]{berthelot1} constructs an isomorphism
between $H^i_{\rig}(X)$ and the cohomology of the complex
\[
\Omega^\cdot_{A/\gotho_K} \otimes_A A^\dagger \otimes_{\gotho_K} K;
\]
the latter is the complex computing Monsky-Washnitzer's ``formal
cohomology'' \cite{vdp}. Moreover, if $F: A^\dagger \to A^\dagger$ is any
ring homomorphism extending $\sigma$ on $\gotho_K$ and lifting the absolute
Frobenius, then pullback by $F$ induces the Frobenius action on rigid cohomology.
\end{defn}

\begin{remark} \label{R:rigid pn}
For example, one can compute from the above description that
the $q$-power Frobenius action on $H^{2i}_{\rig}(\PP^n_{\FF_q})$
consists of multiplication by $q^i$. This also follows from the fact that
the cohomology of projective space is generated by cycle clasess.
\end{remark}

\section{The case of  smooth hypersurfaces}

From the first part of this paper, we obtain
a procedure for bounding from above 
the Picard number and the geometric Picard number of a smooth proper
variety $X$ over a finite field $\FF_q$: compute an approximation modulo
$p^m$, for some $m$, to the matrix via which Frobenius acts on
the rigid cohomology space $H^2(X)$ over $\QQ_q$ (i.e., a ``higher
Cartier matrix''), then use
Algorithm~\ref{algo:picard} to bound the right-hand side of
\eqref{eq:tate3} or \eqref{eq:tate4}, respectively. 
What remains to be
done, for any particular class of varieties,
is to describe how to compute the approximate Frobenius matrix for varieties
in that class, by realizing the constructions described in the second part of
the paper. Here, we describe one such procedure for smooth
hypersurfaces in a projective space, based on work of Griffiths
\cite{griffiths}, and give a few details of an implementation of
this procedure which we have constructed. We also mention some
alternate approaches that we have not (yet) experimented with.

It is  worth noting that much of what is described below generalizes
relatively easily to smooth hypersurfaces, or even smooth complete
intersections, in toric varieties;
for example, Johan de Jong is currently developing an implementation for
hypersurfaces in weighted projective spaces \cite{dejong}.
We have restricted to hypersurfaces
in projective space, and ultimately to surfaces in $\PP^3$, both
for simplicity
of exposition and because that is all that we have attempted to implement
ourselves.

\subsection{A calculation on projective space}


We pause for a brief excursion into the cohomology of sheaves of differentials
on projective space.

\begin{lemma} \label{L:diff exact}
Let $S$ be a scheme. For any positive integer $n$,
there is an exact sequence
of sheaves on $X = \PP^n_S$:
\[
0 \to \Omega^1_{X/S} \to \calO_X(-1)^{n+1} \to \calO_X \to 0.
\]
\end{lemma}
\begin{proof}
A standard calculation: see \cite[Theorem~II.8.13]{hartshorne}.
\end{proof}

The following is due to Bott \cite{bott} over $\CC$ (and can be
deduced over any field of characteristic zero); we were unable to find
a reference for the general version, so we include the easy calculation.
One might like to
think of it as a special case of the Kodaira-Nakano vanishing
theorem, but the latter is not valid over a general base;
Raynaud \cite{raynaud} exhibited a counterexample over a field of
positive characteristic.
\begin{prop} \label{P:bott}
Let $S = \Spec A$ be an affine scheme. For $n$ a
positive integer, put $X = \PP^n_S$. Let $k,p,q$ be integers with
$0 \leq p,q \leq n$. Then $H^q(X, \Omega^p_{X/S}(k)) = 0$ unless
one of the following conditions holds:
\begin{enumerate}
\item[(a)]
$p=q$ and $k=0$, in which case $H^p(X, \Omega^p_{X/S}(k))$ is free of rank $1$
over $A$;
\item[(b)]
$q=0$ and $k>p$;
\item[(c)]
$q=n$ and $k <p-n$.
\end{enumerate}
\end{prop}
\begin{proof}
We will use without comment 
Serre's calculation of the cohomology of $\calO_X(k)$, in the form
of the statement of the lemma in all cases where $p=0$; for this, see
\cite[Theorem~III.5.1]{hartshorne} (and suppress the superfluous
noetherian hypothesis) or
\cite[Proposition~2.1.12]{ega3-1}.

Take the $p$-th exterior power of the exact sequence in
Lemma~\ref{L:diff exact}, then twist by $k$, to obtain
\[
0 \to \Omega^p_{X/S}(k) \to
\left(\wedge^p (\calO_X(-1)^{n+1}) \right) (k) \to 
\Omega^{p-1}_{X/S}(k) \to 0,
\]
which we may rewrite as
\begin{equation} \label{eq:serre exact}
0 \to \Omega^p_{X/S}(k) \to \calO_X(-p+k)^{\binom{n+1}{p}}
\to \Omega^{p-1}_{X/S}(k) \to 0.
\end{equation}
We now proceed by inspecting part of the long exact sequence
in cohomology of \eqref{eq:serre exact}:
\begin{equation} \label{eq:serre exact2}
H^{q-1}(X, \calO_X(-p+k))^{\binom{n+1}{p}}
\to H^{q-1}(X, \Omega^{p-1}_{X/S}(k)) 
\to H^{q}(X, \Omega^p_{X/S}(k)) 
\to
H^{q}(X, \calO_X(-p+k))^{\binom{n+1}{p}}.
\end{equation}
In case $k < p-q$ and $0 \leq  q< n$, then the outside terms of
\eqref{eq:serre exact2} vanish, so
$H^{q-1}(X, \Omega^{p-1}_{X/S}(k))$ and $H^{q}(X, \Omega^p_{X/S}(k))$ are
isomorphic. We obtain the vanishing of $H^{q}(X, \Omega^p_{X/S}(k))$ in
this case by successively decreasing both $p$ and $q$ until the step
when at least
one of them 
goes negative, at which moment the correct formal interpretation
of \eqref{eq:serre exact2} yields the desired vanishing.
Similarly, in case
$k > p -q$ and $0 < q \leq n$, we obtain vanishing of
$H^q(X, \Omega^p_{X/S}(k))$ by successively increasing both $p$ and $q$
until the step when at least one of them exceeds $n$.

In case $k = p-q$,
the outside terms of \eqref{eq:serre exact2} still vanish as long as
$1 \leq q \leq n$. If $k \neq 0$, then successively decreasing $p$ and $q$
still eventually manages to send one of them below zero, so we get
the desired vanishing.
If $k=0$, we instead hit $H^0(X, \calO_X)$ which is free of rank 1.
This yields all of the desired results.
\end{proof}
\begin{cor} \label{C:middle cohom}
Let $S = \Spec A$ be an affine scheme, let $n$ be a positive integer, and
put $X = \PP^n_S$.
Let $k$ be a nonnegative integer, let $Z$ be a smooth hypersurface in $X$,
and define the complex
$C^p = \Omega^p_{X/S}((k+p)Z)$ with the evident differentials $d$. 
Then the hypercohomology of $C^\cdot$
coincides with the homology of the complex $H^0(X, C^\cdot)$. In particular,
\[
\HH^n(X, C^\cdot) = H^0(X, C^n) / dH^0(X, C^{n-1}).
\]
\end{cor}
\begin{proof}
We compute $\HH^n(X, C^\cdot)$ using a spectral sequence with
$E_1^{pq} = H^q(X, C^p)$. By Proposition~\ref{P:bott}
and the fact that $C^p \cong \Omega^p_{X/S}((k+p)\deg(Z))$,
we have $E_1^{pq} = 0$
for $q>0$. Hence the sequence degenerates at $E_2$ and yields the claim.
\end{proof}
\begin{cor} \label{C:pn}
For any affine scheme $S = \Spec A$,
\[
H^i_{\dR}(\PP^n_S) \cong \begin{cases} A & i = 0, 2, \dots, 2n \\
0 & \mbox{otherwise}.
			 \end{cases}
\]
\end{cor}
\begin{proof}
By Proposition~\ref{P:bott},
$H^p(\PP^n_S, \Omega^q_{\PP^n_S/S})$ is free of rank 1 over $A$ if
 $p = q \in \{0, \dots, n\}$ and is otherwise zero. Hence the
spectral sequence computing hypercohomology degenerates at $E_1$,
yielding the desired result.
\end{proof}

\subsection{Cohomology of smooth hypersurfaces (after Griffiths)}

The middle cohomology of a smooth hypersurface in a projective space was
described by Griffiths \cite{griffiths} using mostly analytic arguments (i.e.,
working over $\CC$ and invoking GAGA).
One can reconstruct these results algebraically; we will not do so explicitly,
but we will use algebraic techniques later to extract arithmetic information.

\begin{notation}
Throughout this section,
let $K$ be a field of characteristic zero, put $X = \PP^n_K$, let $Z$ be a 
smooth hypersurface of degree $d$
in $X$, defined by the homogeneous polynomial $P \in K[x_0, \dots, x_n]$,
and put $U = X \setminus Z$.
By the Lefschetz hyperplane theorem, the map $H^i_{\dR}(X/K) \to H^i_{\dR}(Z/K)$ 
induced by the
inclusion $Z \hookrightarrow X$ is bijective for $i \leq n-2$ and
injective for $i=n-1$. By Corollary~\ref{C:pn} and
Poincar\'e duality, we thus conclude that
for $0 \leq i \leq 2n-2$ with $i \neq n-1$, we have
\[
\dim_K H^i_{\dR}(Z/K) = \begin{cases} 1 & \mbox{$i$ even} \\
0 & \mbox{$i$ odd}.
		      \end{cases}
\]
In particular, the only cohomology group of $Z$ requiring further inspection
is the middle cohomology $H^{n-1}_{\dR}(Z/K)$.
\end{notation}

The following result lets us shift attention from $Z$ to the ambient
projective space $X$, where it is easier to make calculations.
Remember that $H^n_{\dR}((X,Z)/K) \cong H^n_{\dR}(U/K)$ by
Corollary~\ref{C:deligne}.

\begin{prop} \label{P:induced exact}
The exact sequence of Proposition~\ref{P:excision} induces exact sequences as 
follows: if $n$ is even, then
\[
0 \to H^{n}_{\dR}(U/K) \to H^{n-1}_{\dR}(Z/K) \to 0
\]
is exact; if $n$ is odd, then
\[
0 \to H^{n}_{\dR}(U/K) \to H^{n-1}_{\dR}(Z/K) \to H^{n+1}_{\dR}(X/K) \to 0
\]
is exact.
\end{prop}
\begin{proof}
See \cite[(10.16)]{griffiths}.
\end{proof}

\begin{defn} \label{D:reduction}
Put
\[
\Omega = \sum_{i=0}^n (-1)^i x_i \,dx_0 \wedge \cdots \wedge
\widehat{dx_i} \wedge \cdots \wedge dx_n,
\]
where the hat denotes omission. 
One then calculates as in \cite[\S 4]{griffiths} that
$H^n_{\dR}((X,Z)/K)$ is isomorphic to the quotient of the group of
$n$-forms $A \Omega/P^{m}$, where $m$ is a positive integer and $A
\in K[x_0, \dots, x_n]$ is homogeneous of degree $md - n - 1$, by the subgroup
generated by
\[
\frac{(\del_i A) \Omega}{P^m} - m \frac{A (\del_i P) \Omega}{P^{m+1}}
\]
for all nonnegative integers $m$, all $i \in \{0, \dots, n+1\}$, and all
homogeneous polynomials $A \in K[x_0, \dots, x_n]$ of degree
$md - n$. (Here $\del_i$ is shorthand for $\frac{\partial}{\partial x_i}$.)
\end{defn}

\begin{remark}
Note that $H^n_{\dR}(U/K)$ admits a natural filtration whose $i$-th
step consists of those classes represented by forms $A\Omega/ P^{m}$
for some integer $m \leq i+1$. In fact, this filtration is 
induced by the Hodge filtration on $H^{n-1}_{\dR}(Z/K)$
\cite[\S 10]{griffiths}.
\end{remark}

\begin{remark} \label{R:reduction}
The description in Definition~\ref{D:reduction}
gives rise to a natural ``reduction of poles''
procedure for computing in $H^n_{\dR}(U)$, sometimes referred to as the
Griffiths-Dwork method. First, one writes down
a basis: for $h=1, \dots, n$, one finds monomials of degree
$hd-n-1$ which generate the quotient of the space of all such monomials
by the multiples of $\del_0 P, \dots, \del_n P$. Then, to write the class of a
given form $A\Omega/P^m$ in terms of these, one uses a Gr\"obner basis
division procedure to write $A$ as a linear combination of
$\del_0 P, \dots, \del_n P$ (plus basis elements if $m \leq n$), 
then reduces the pole order. The fact that it is always possible to perform
this reduction follows
from a theorem of Macaulay \cite[\S 4]{griffiths} or from a sheaf-theoretic
reinterpretation \cite[\S 10]{griffiths}.
\end{remark}

\subsection{The $p$-adic cohomology interpretation}

We now use the previous subsection to describe the $p$-adic cohomology of a smooth
hypersurface in $\PP^n_{\FF_q}$. Note that while the $p$-power Frobenius
action on rigid cohomology of a variety over $\FF_q$ will only be semilinear,
the $q$-power Frobenius action will be linear over $\QQ_q$.

We start by setting notation for the rest of the section.
\begin{notation} \label{N:hypersurface}
Let $Z$ be a smooth hypersurface of degree $d$ 
in $X = \PP^n_{\FF_q}$, for
$\FF_q$ a finite field of characteristic $p>0$, defined by the 
homogeneous polynomial $P(x_0, \dots, x_n) \in
\FF_q[x_0, \dots, x_n]$. Choose a lift $\frakP(x_0, \dots, x_n) \in
\ZZ_q[x_0, \dots, x_n]$ of $P$ to a homogeneous polynomial of the
same degree $d$, let $\frakZ$ be the zero locus of $\frakP$ in
$\frakX = \PP^n_{\ZZ_q}$, and put $\frakU = \frakX \setminus \frakZ$.
Put $\tilde{X} = \frakX_{\QQ_q}$, $\tilde{Z} = \frakZ_{\QQ_q}$,
and $\tilde{U} = \frakU_{\QQ_q}$; also write $\tilde{P}$ for
$\frakP$ as an element of $\QQ_q[x_0, \dots, x_n]$.
\end{notation}

By Definition~\ref{D:rigid}, $H^i_{\rig}(Z) \cong 
H^i_{\crys}(Z) \otimes_{\gotho_K} K$;
adding Proposition~\ref{P:liftable}(a), we get $H^i_{\rig}(Z) \cong
H^i_{\dR}(\frakZ/\ZZ_q) \otimes_{\ZZ_q} \QQ_q \cong 
H^i_{\dR}(\tilde{Z}/\QQ_q)$.
Adding the Lefschetz hyperplane theorem, we obtain
$H^i_{\rig}(Z) \cong H^i_{\rig}(X)$ for
$i \leq n-2$. 
Adding Poincar\'e duality, we see that for $0 \leq i \leq 2n-2$
with $i \neq n-1$, if $i$ is odd, then $H^i_{\rig}(Z) = 0$, while if 
$i$ is even, then $H^i_{\rig}(Z)$ is one-dimensional
and the $q$-power Frobenius acts by multiplication by $q^{i/2}$.
Moreover,
by Proposition~\ref{P:induced exact} 
and Proposition~\ref{P:liftable}(b),
we have Frobenius-equivariant exact sequences
\[
0 \to H^n_{\rig}(U) \to H^{n-1}_{\rig}(Z)(-1) \to 0
\]
if $n$ is even and
\[
0 \to H^n_{\rig}(U) \to H^{n-1}_{\rig}(Z)(-1) \to H^{n+1}_{\rig}(X) \to 0
\]
if $n$ is odd. That is, 
$H^n_{\rig}(U)(1)$ coincides with $H^{n-1}_{\rig}(Z)$ except that
if $n$ is odd, its generalized eigenspace for Frobenius
of eigenvalue $q^{(n-1)/2}$ has dimension one less. (That is,
$H^n_{\rig}(U)$ is the \emph{primitive middle cohomology} of $Z$, 
i.e., the part not explained by the Lefschetz hyperplane theorem.)

Thanks to Berthelot's comparison theorems (see Definition~\ref{D:monsky}),
we can describe the Frobenius action on $H^n_{\rig}(U)$ as follows.

\begin{defn}
Let $v$ denote the Gauss valuation on polynomials over $\QQ_q$; that is,
$v(\sum c_I x^I) = \min_I \{v_p(c_I)\}$, where $v_p$ denotes the
$p$-adic valuation on $\QQ_q$ normalized by $v_p(p) = 1$.
\end{defn}

\begin{defn}
Let $R$ denote the ring of formal sums $\sum_{i=0}^\infty A_i \tilde{P}^{-i}$,
where $A_i \in \QQ_q[x_0,\dots,x_n]$ is homogeneous of degree $di$, and
\[
\liminf_{i \to \infty} \frac{v(A_i)}{i} > 0.
\]
That is, the valuations of the $A_i$ grow at least 
linearly in $i$. Define the ring
map $F: R \to R$ by formally setting $F(x_i) = x_i^q$ for $i=0,\dots, n$,
and
\[
F(\tilde{P}^{-1}) = \tilde{P}^{-q} \left( 1 + p \frac{F(\tilde{P}) - \tilde{P}^q}{p \tilde{P}^q} \right)^{-1},
\]
expanding the parenthesized expression as a binomial series.
Formally extend $F$ to $n$-forms by setting
\[
F(A \Omega) = F(A x_0 \cdots x_n) F(x_0^{-1} \cdots x_n^{-1} \Omega),
\]
where $F(dx_0/x_0) = q\,dx_0/x_0$ and so forth.
As noted in Definition~\ref{D:monsky}, this ring map
induces the $q$-power Frobenius in rigid cohomology on $H^n_{\rig}(U)
\cong H^n_{\dR}((\tilde{X},\tilde{Z})/\QQ_q) \cong H^n_{\dR}(\tilde{U}/\QQ_q)$.
\end{defn}

\subsection{Precision estimates}

The plan now is to compute an approximation to the Frobenius action on
$H^n_{\rig}(U)$ by applying a truncation of $F$ to a basis of 
$H^n_{\dR}(\tilde{U})$ and using the ``reduction of poles'' process
(Remark~\ref{R:reduction}) to put the results back in terms of the basis.
To do this, we must produce \emph{effective} bounds on the amount of
$p$-adic precision needed to keep the error introduced by the truncation
to a particular amount. One can derive general bounds easily from the theory of
$p$-adic cohomology, but for effective bounds we must work a bit harder.
This is analogous to the precision analysis in \cite{me-count}, but the
higher-dimensional situation we are in makes things a bit more technical.

Our first order of business is to relate a basis, which is natural from
the point of view of reduction of poles, 
to the integral structure on de Rham cohomology. 

\begin{defn}
Let $H$ be the image of $H^n_{\dR}((\frakX, \frakZ)/\ZZ_q)$
in $H^n_{\dR}((\tilde{X},\tilde{Z})/\QQ_q)$; we refer to elements of $H$
as \emph{integral} elements of 
$H^n_{\dR}((\tilde{X},\tilde{Z})/\QQ_q)$, or of 
$H^n_{\dR}(\tilde{U}/\QQ_q)$,
or of $H^n_{\rig}(U)$.
\end{defn}

\begin{defn} \label{D:basis}
Let $B$ denote a basis of $H^n_{\dR}(\tilde{U}/\QQ_q)$ obtained as follows.
For $h=1,\dots,n$, form the quotient of the space of homogeneous polynomials
in $\FF_q[x_0,\dots,x_n]$ of degree $hd-n-1$ by the multiples of
$\del_0 P, \dots, \del_n P$. Find monomials in $\FF_q[x_0,\dots,x_n]$
which project onto a basis of this quotient, then lift these
monomials to monomials in $\ZZ_q[x_0,\dots,x_n]$.
For each such lift $\mu$, include $\mu \Omega/\tilde{P}^h$
in $B$. Let $V \subset H^n_{\dR}(\tilde{U}/\QQ_q)$ denote
 the $\ZZ_q$-span of $B$. (This is not the only logical 
choice; see Remark~\ref{R:alt basis}.)
\end{defn}

\begin{lemma} \label{L:compare basis}
Let $W$ be the $\ZZ_q$-span in
$H^n_{\dR}(\tilde{U}/\QQ_q)$
of the elements $\mu \Omega/\tilde{P}^h$, for $h \in \{1,\dots,n\}$
and $\mu \in
\ZZ_q[x_0,\dots,x_n]$ homogeneous of degree $hd-n-1$.
Then
\[
H \subseteq W.
\]
\end{lemma}
\begin{proof}
Map the complex $\Omega^{\cdot}_{(\frakX,\frakZ)/\ZZ_q}$
into the complex $C^\cdot$ with
\[
C^i = \Omega^{i}_{\frakX/\ZZ_q}(i\frakZ),
\]
then invoke Corollary~\ref{C:middle cohom}.
\end{proof}

\begin{remark}
In general, we have $H \subseteq W$ by Lemma~\ref{L:compare
basis}, $V \subseteq W$ evidently, 
and $(n-1)!W \subseteq V$ by inspection of the
reduction process. 
In the special case $p \geq n$, though,
it will follow from Corollary~\ref{C:compare basis2}
that $H = V = W$. 
In this case, we also know that
$H^{n-1}_{\dR}(\frakZ/\ZZ_q)$ is torsion-free because
of the degeneration of the Hodge-de Rham spectral sequence. See
\cite{illusie} for further discussion.
\end{remark}

We next consider the loss of $p$-adic precision incurred when one
reduces a given differential into standard form.

\begin{defn}
For $m$ a positive integer, let $f(m)$ be the smallest integer $t$
with the following property: for each 
form $\omega = A \Omega/\tilde{P}^m$ with $A \in \ZZ_q[x_0,\dots,x_n]$
homogeneous of degree $md-n-1$, 
$p^t \omega$ represents an element of $H$. The following relations are evident
(using the relations in cohomology):
\begin{align*}
f(1) &= 0 \\
f(m) &\leq f(m+1) \\
f(m) &= f(p \lceil m/p \rceil) \qquad (m \geq n).
\end{align*}
\end{defn}

\begin{prop} \label{P:prec bound1}
For each $m > 0$,
\[
f(m) \leq \sum_{i=1}^n \lfloor \log_p \max\{1,m-i\} \rfloor.
\]
\end{prop}
\begin{proof}
Apply Theorem~\ref{thm:deligne};
we get the $n$ terms of the sum by feeding the
cohomology sheaves into the spectral sequence computing hypercohomology,
keeping in mind Remark~\ref{R:lower pole order}.
(We don't get an $(n+1)$-st term because the map of Theorem~\ref{thm:deligne}
on zero-th cohomology sheaves is an isomorphism, so those do not contribute.)
\end{proof}
\begin{cor} \label{C:compare basis2}
We have $f(m) \leq v_p((m-1)!)$ for all $m \geq 0$. In particular,
in the notation of Lemma~\ref{L:compare basis}, we have 
$(n-1)! W \subseteq H \subseteq W$.
\end{cor}

This bound is asymptotically $n \log_p(m)$, which for our application to
surfaces will be a bit too large to be practical. Fortunately we can 
shave it down a bit.

\begin{defn}
For $m,i$ integers with $i \geq 0$, let $g(m,i)$ be the $p$-adic valuation of
$\binom{-m}{i}$. By a standard argument attributed to Kummer, $g(m,i)$
equals the number of carries in the base $p$ addition of $i$ and $-m-i$.
\end{defn}

\begin{prop} \label{P:prec bound2}
Let $m$ be a positive integer.
Put
\[
N = \max_{\ell>0} \{ 
f((m+\ell)p) - \ell - g(m,\ell) \}.
\]
Then
\[
f(mp) \leq \max\{N, n-1+f(m)\}.
\]
\end{prop}
\begin{proof}
Given a form $A \Omega/\tilde{P}^{mp}$ with $A \in \ZZ_q[x_0,\dots,x_n]$, 
separate the monomials of $A$
depending on the reductions modulo $p$ of their exponents.
Let $B$ be the sum of the monomials whose exponents are all congruent to
$p-1$ modulo $p$, and put $C = A - B$. Then
$C \Omega/\tilde{P}^{mp}$ is cohomologous to a form
$mp D \Omega/\tilde{P}^{mp+1}$ with $D \in \ZZ_q[x_0, \dots,x_n]$,
because
\[
\frac{x_i^j \Omega}{\tilde{P}^{mp}} \equiv \frac{mp}{j+1} \frac{x_i^{j+1}
(\del_i \tilde{P}) \Omega}{\tilde{P}^{mp+1}}
\]
in cohomology. Such a form reduces to an element of $p^{-t} H$ 
for $t = f(mp+1) - 1 - g(m,1) \leq N$.
On the other hand, we can write $B = p^{-n}F(D) + E$
with $D$ having pole order $m$ and
$E = 
\sum_{\ell>0} p^\ell \binom{-m}{\ell} E_\ell \Omega/\tilde{P}^{(m+\ell)p}$
for some $E_{\ell} \in \ZZ_q[x_0,\dots,x_n]$.
Since $p^{-1} F$ acts on $H$ (by comparison with
$H^{n-1}_{\dR}(\frakZ/\ZZ_q)$ via Proposition~\ref{P:liftable}),
$D$ reduces to an element of $p^{-t} H$ for $t = n-1+f(m)$,
while the $\ell$-th summand of $E$ reduces to an element of $p^{-t} H$ for
$t = f((m+\ell)p) - \ell - g(m,\ell) \leq N$.
This yields the claim.
\end{proof}

Proposition~\ref{P:prec bound2} plus any sublinear bound on $f$
suffices to give an upper bound of the form $(n-1) \log_p(m)$ plus 
a constant. However, for implementation purposes, it is important
to control that additive constant as much as possible; this can be done
using an iterative algorithm as follows.

\begin{algo} \label{algo:prec}
Define the function
\[
f_0(m) = \sum_{i=1}^n\lfloor \log_p 
\max\{1,m-i\} \rfloor.
\]
Given an input positive integer $m$ and an input array $A_0$ such that
$f(j) \leq A_0(j)$ for each $j$ for which $A_0(j)$ exists, 
the following algorithm, 
if it terminates, returns an array $A$ of length at least $m$,
such that $f(j) \leq A(j)$ for each $j$ for which $A(j)$ exists,
and $A(j) \leq A_0(j)$ for each $j$ for which $A(j)$ and $A_0(j)$ both exist.
\begin{enumerate}
\item
Create an array 
\[
A(j) = \begin{cases} A_0(j) & \mbox{$A_0(j)$ exists} \\
f_0(j) & \mbox{otherwise} \end{cases}
\qquad (j=1, \dots, \min\{n,m\}).
\]
\item
Make a copy $A'$ of $A$. Put $j=n+1$.
\item
If $A(j)$ is not defined, check whether $A$ and $A'$ are identical
(of the same length). If so, return $A$ and STOP. Otherwise,
go to step 2. (If $A(j)$ is defined, continue to step 4.)
\item
Put $j_1 = p \lceil j/p \rceil$,
$N = n - 1 + A(j_1/p)$, and $\ell = 1$.
\item
If $np < j_1+\ell p$ and 
$n \log_p (j_1 + \ell p) - \ell \leq N$, then set 
\[
A(j)=A(j+1) = \cdots= A(j_1) = \min\{N, f_0(j_1)\},
\]
replace $j$ by $j_1 + 1$, and go to step 3.
\item
If $f_0(j_1 + \ell p) -\ell - g(j_1,\ell) \leq N$,
then replace $\ell$ by $\ell+1$ and go to 
step 5.
\item
Extend $A$ using the formula $A(i)= f_0(i)$ if needed to ensure that
$A(j_1 + \ell p)$ exists. 
Replace $N$ by $\max\{A(j_1 + \ell p)-\ell-g(j_1,\ell),N\}$, replace $\ell$ by
$\ell+1$, and go to step 5.
\end{enumerate}
\end{algo}
\begin{proof}
What we show is that at every stage, whenever some $A(j)$ is defined,
we have $f(j) \leq A(j)$. This holds whenever an $A(j)$ is instantiated
by Proposition~\ref{P:prec bound1}.
In step 5, the quantity $n \log_p (j_1 + \ell p) - \ell$ is
a decreasing function of $\ell$ for $np< j_1 + \ell p$; it is also an
upper bound for $f(j_1+\ell p)- \ell - g(j_1,\ell)$ by 
Proposition~\ref{P:prec bound1} again. 
The property that $f(j) \leq A(j)$ is preserved in step 5 thanks to
Proposition~\ref{P:prec bound2}.
\end{proof}

\begin{remark}
One can prove termination of Algorithm~\ref{algo:prec} and control its runtime
with a bit of effort, but in practice it suffices to simply let it run
until either the process terminates, or one goes through a set number
of iterations of step 3 (we used 20 iterations in our examples).
\end{remark}

\begin{remark} \label{R:alt basis}
If in Definition~\ref{D:basis} we had taken $B$ to consist of the
elements $(h-1)!\mu \Omega/\tilde{P}^h$ instead of $\mu \Omega/\tilde{P}^h$,
we would get a $\ZZ_q$-span $V'$ satisfying $(n-1)! W \subseteq V' \subseteq H
\subseteq W$. However, in practice this choice appears to give slightly
less control of the precision loss.
\end{remark}

\subsection{Summary of the algorithm}

To summarize, we describe how to assemble an algorithm for computing
an approximate Frobenius matrix on $H^n_{\rig}(U)(-1)$.

To start with, write down a basis $B$ of rigid cohomology as in
Definition~\ref{D:basis}. Say we want to compute the Frobenius matrix
on this basis modulo $p^r$. Choose the integer $s$ using
the following algorithm:
\begin{enumerate}
\item
Put $s=r$. Create an empty array $A$.
\item
Put $j=s-n+1$.
\item
If $j>0$ and $n \log_p(p(n+j)-1) \leq n-1+j-r$, then
return $s$ and STOP.
\item
Replace $A$ by the result of Algorithm~\ref{algo:prec} applied with
inputs $p(n+j)$ and $A$.
If $A(p(n+j)) > n-1+j-r$, then
replace $s$ by $s+1$ and go to step 2. Otherwise, replace $j$ by $j+1$
and go to step 3.
\end{enumerate}

For each basis element, compute the image of Frobenius truncating all
terms which vanish modulo $p^s$, then use reduction of poles
(Remark~\ref{R:reduction}; see also Remark~\ref{R:reduction2} below) 
to express the result in terms of the basis elements.
The verification that this gives enough precision is straightforward:
the Frobenius image
of a monomial $\mu \Omega/\tilde{P}^h$ has the form
$\sum_{j=0}^\infty p^{n+j-1} B_j \Omega/\tilde{P}^{p(h+j)}$,
so
we need to ensure that 
\[
n-1+j-f(p(h+j)) \geq r \qquad (h \leq n, j \geq s-n+1),
\]
but checking for $h=n$ implies the same for $h \leq n$,
and the upper bound from Proposition~\ref{P:prec bound1} allows us to
truncate in step 3.

\begin{remark} \label{R:parallel}
This algorithm readily admits some parallelization, as one can
compute the reductions of the Frobenius images of different basis elements
independently. The experimental results we describe later do not depend
on this capability, but it may be useful for larger examples.
\end{remark}

\begin{remark} \label{R:reduction2}
There are at least four distinct ways to carry out the reduction of poles
implied by Remark~\ref{R:reduction}.
\begin{itemize}
\item
One can simply perform the entire calculation over $\QQ$, then interpret
the final result modulo the appropriate power of $p$. This was our
first choice, implemented in \textsc{Singular}, 
but it leads to undesirable intermediate coefficient blowup.
\item
One can perform the calculation over $\ZZ/p^s \ZZ$ by using integral
analogues of Gr\"obner basis arithmetic as implemented in \textsc{Magma}: 
to do this, one must postpone
the division by $m-1$ implied in reducing the pole order from $m$ to $m-1$
until the end of the calculation. One must also remember to use the trivial
relation $A\Omega/\tilde{P}^m = A\tilde{P}\Omega/\tilde{P}^{m+1}$, as 
$P$ may not be generated by its partial derivatives (in case $p$ divides
$\deg(P)$). This was our second choice, and is used in our current
implementation.
\item
One can perform the calculation over $\FF_{\ell}$ for any prime $\ell$ 
of good reduction of the lifted hypersurface. This would allow the use of
the more efficient polynomial arithmetic of \textsc{Singular} over 
\textsc{Magma}; however, in order to recover the desired answer (by working
modulo many small $\ell$ and using the Chinese remainder theorem) 
one would need a bound on the heights of the entries of the resulting matrix.
We have not attempted this method.
\item
One can perform the calculation over $\CC$ by 
going through the Betti-de Rham
comparison as in \cite{griffiths}: numerically integrate each truncated 
Frobenius
image against a basis of homology, then perform a lattice reduction to express
these in terms of the integrals of a basis of cohomology. 
Again, this requires height bounds, and again we have not attempted this method.
\end{itemize}
The last two methods share ideas
with the method used by Edixhoven et al \cite{edixhoven}
to compute coefficients of the $\Delta$ modular form
(and by extension with the Schoof-Elkies-Atkin method for computing zeta 
functions of elliptic curves).
\end{remark}

\subsection{Alternate algorithms}
\label{subsec:alternate}

As the experimental data in the final part of the paper suggests,
computing approximate Frobenius matrices in
$p$-adic cohomology in 
the manner we have suggested is rather laborious. There are several other
ways one might perform this computation; we mention some of these in passing,
noting that any of them can be used together with the first part of the paper
to give a test for low Picard number. (The arguments in the second part of
the paper, notably Theorem~\ref{thm:deligne}, may help in the analysis
of precision loss in such algorithms.)

The shift from a hypersurface to its affine complement amounts to an increase
by one in the dimensions of the varieties under consideration, and in the
number of variables in the polynomial rings in which the calculations take 
place. (In reality there is one more variable even than that, but this is 
merely because we are working with homogeneous polynomials.)
That shift turns out to be costly, so one would ideally like to avoid
it. This appears to be possible for so-called nondegenerate hypersurfaces,
those which together with the coordinate hyperplanes and the hyperplane at
infinity form a normal crossings divisor. For curves, this has been proposed by
Castryck, Denef, and Vercauteren \cite{cdv}, 
but the method extends readily
to higher dimensions. We made a cursory
attempt to implement this method for surfaces in $\PP^3$,
but our results were inconclusive: the
additional complexity in the method (especially in lifting Frobenius,
which would now be done on an affine piece of the original hypersurface
rather than on an affine complement) seemed to introduce 
large constant factors that interfered with the
asymptotic improvements at the scale at which our calculations took place.
Nonetheless, we think the method deserves further study.

A better approach may be to use d\'evissage: write the given
 surface as a fibration
of curves, compute the higher direct images of the constant sheaf, then
compute the cohomology of these. This has been suggested by Lauder, who
has implemented this in some examples with good results \cite{lauder-recur}.

Yet another approach is to avoid directly computing the cohomology of
the particular hypersurface of interest, by instead putting it into a
pencil with one member chosen to be smooth with extra automorphisms. One
can then compute its Frobenius matrix more easily, then use that data as
the ``initial condition'' in the differential equation relating the Gauss-Manin
connection of the pencil with its global Frobenius action.
This is the ``deformation'' method of Lauder \cite{lauder-def,
lauder-def2}; it has been tested experimentally for families
of elliptic curves by Gerkmann  \cite{gerkmann2}
and has been theoretically analyzed for
hyperelliptic curves by Hubrechts \cite{hubrechts} (where it already
gives some improvement over the direct method), but we are not aware
of any work in higher dimensions besides Lauder's original papers.

\section{Implementation details}

In this section, we describe an implementation that implements a special case
of the algorithm we have described, and give some experimental results.
One glaring omission is that we do not make a complexity analysis; this is only
partly out of laziness. The other reason is that Gr\"obner basis calculations
in general have extremely bad worst-case performance; we are not in the worst
case here, but we would have to look closely at what we are using to obtain
complexity estimates. Since the purpose of this paper is instead
to demonstrate
the practicality of these methods, we do not carry out such intricate analysis
here.

\subsection{Implementation notes}

Using the \textsc{Magma} algebra
system \cite{magma},
we have developed an implementation of the
methods of this paper, to obtain an algorithm for computing approximate
Frobenius matrices in rigid cohomology
for smooth surfaces in projective 3-space over a prime finite field.
See \cite{akr} for the source code.

We have tested this implementation on the computer \texttt{dwork.mit.edu},
a Sun workstation with dual Opteron 246 CPUs running at 2 GHz, with
access to 2GB of RAM. Although
these CPUs are 64-bit processors, these experiments were conducted in 32-bit
mode under Red Hat Enterprise Linux 4. 
Each individual surface was run on a single CPU with no use of parallelism
(see Remark~\ref{R:parallel}), and timings
are reported in CPU seconds; memory usages are reported in megabytes.
Beware that these should only be taken as order-of-magnitude indications:
there are slight variations from run to run of a single example, and there
are much bigger variations within classes of examples, probably
arising from the vagaries of Gr\"obner basis arithmetic.

Some initial experiments were also conducted using the \textsc{Singular}
algebra system \cite{singular}. The main downside with using \textsc{Singular}
for this calculation is that it only treats polynomials over a field; while
one obtains correct answers by reducing poles over $\QQ$ and then reducing
modulo a power of $p$, the resulting
calculations experience unacceptable intermediate coefficient blowup. By
contrast, the \textsc{Magma} implementation uses a Gr\"obner basis 
implementation over $\ZZ/p^m\ZZ$, which avoids the coefficient blowup.
(Compare Remark~\ref{R:reduction2}.)

In the subsequent sections, we describe some examples computed using
this implementation. These examples were chosen to be ``generic'', without
special geometric properties: their coefficients were chosen at random with
a bias towards zero coefficients. The bias somewhat simplifies the
Gr\"obner basis calculations. 

It is worth noting that one can use a
``prescreening'' strategy to find such examples: deliberately compute
approximations with \emph{not enough} initial precision, then revisit 
the examples that appear to work with a provably sufficient amount of 
precision. We suspect that this works because
our precision estimates do not give
a complete picture of how quickly the $p$-adic approximations are
converging; see Remark~\ref{R:extra prec} for an instance of this.

\subsection{Example: degree 4 over $\FF_3$}

We start with a careful analysis of an example in what is possibly
the simplest nontrivial case. Namely, surfaces of degree $1$ and $2$
over a finite field $\FF_p$ are isomorphic to $\PP^2$ and $\PP^1 \times \PP^1$,
whose zeta functions are known,
while surfaces of degree 3 have all cohomology generated by the classes
of the 27 lines on the surface over $\overline{\FF_p}$, so the zeta function
can be computed from the Galois action on these lines.
(Note that while a smooth cubic surface over $\FF_p$ has geometric
Picard number 7, its arithmetic Picard number can equal 1;
see \cite{zarzar} for an example.)

A smooth surface in $\PP^3$ of degree 4 is a K3 surface, whose middle
cohomology has dimension 22 and Hodge numbers 1, 20, 1. (Remember that
we compute using \emph{primitive} middle cohomology, in which the dimension
and the central Hodge number are both decreased by 1.) We will exhibit
an example of provable arithmetic Picard number 1 over $\FF_3$; the more
natural first choice $\FF_2$ actually turns out to be somewhat trickier
(see Section~\ref{subsec:f2}).

\begin{example} \label{exa:f3 quartic}
The smooth quartic surface over $\FF_3$ defined by the polynomial
\[
x^4 - xy^3 + xy^2w + xyzw + xyw^2 -
    xzw^2 + y^4 + y^3w - y^2zw + z^4 + w^4
\]
was found to have Picard number 1 by computing a Frobenius matrix
modulo $3^3$. 
To carry out this calculation provably using the optimal bound extracted
from Algorithm~\ref{algo:prec}, it was necessary to truncate differentials
modulo $3^{7}$; using only Proposition~\ref{P:prec bound1} would have
required working modulo $3^{12}$. 
\end{example}

\begin{remark}
For comparison, Table~\ref{tab:deg4f3} gives some timings and memory usages
for various levels of initial and
final precision in Example~\ref{exa:f3 quartic}. 
When a final precision is specified, that means the
initial precision in that row is the minimum needed to guarantee that
final precision in the Frobenius matrix under Algorithm~\ref{algo:prec};
however, see Remark~\ref{R:extra prec}.
\begin{table}[ht] 
\caption{Timings for a smooth quadric over $\FF_3$.}
\label{tab:deg4f3}
\begin{center}
\begin{tabular}{|c|c|c|c|}
\hline
Final precision & Initial precision & CPU sec & MB \\
\hline
$3^2$ & $3^6$ & 227 & 37 \\
$3^3$ & $3^7$ & 731 & 53 \\
--- & $3^8$ & 907 & 64 \\
--- & $3^9$ & 4705 & 124 \\
$3^4$ & $3^{10}$ & 13844 & 906  \\
$3^5$ & $3^{11}$ & 15040 & 1103\\
$3^6$ & $3^{12}$ & 40144 & 1795 \\
\hline
\end{tabular}
\end{center}
\end{table}
\end{remark}

\begin{remark} \label{R:extra prec}
One reality check on Example~\ref{exa:f3 quartic}
is to verify the implied compatibilities
between the answers computed to various $p$-adic precisions; that is,
the computed Frobenius matrices for any two rows in Table~\ref{tab:deg4f3}
should agree modulo the final precision of the earlier row. This turns
out to be true in a strong fashion: some of the calculations
are even more accurate than predicted. Namely, 
we obtain the correct
Frobenius matrices modulo $3^3, 3^4, 3^5, 3^6$ using initial precisions 
$3^6, 3^7, 3^9, 3^{10}$,
respectively.
\end{remark}

\begin{remark} \label{R:compare charpoly}
Another reality check on Example~\ref{exa:f3 quartic}
is a comparison of initial
coefficients of the predicted zeta functions against those coefficients
determined exactly by actually enumerating rational points.
Over $\FF_{3^i}$ for $i=1,2,3,4,5$, the respective numbers 
of rational points on the
surface are
\[
8, 80, 713, 6836, 58868;
\]
this counting took several hours on a laptop computer using a simple-minded
\textsc{Magma} program.
From the formula for the zeta function (Definition~\ref{D:zeta}), 
and the fact that it has the form $(Q(T)(1-T)(1-3T)(1-9T))^{-1}$ where
$Q(T)$ has degree 21, we deduce that the characteristic polynomial of 
$3^{-1} F$ on primitive middle cohomology begins
\[
\frac{1}{3} (3 T^{21} + 5 T^{20} + 6 T^{19} + 7 T^{18} + 
5 T^{17} + 4 T^{16} + \cdots).
\]
On the other hand, applying Remark~\ref{R:char poly} to the Frobenius
matrix computed with final precision $3^6$, we determine that the
same characteristic polynomial is congruent modulo $3^4$ to
\begin{gather*}
\frac{1}{3} (3 T^{21} + 5 T^{20} + 6 T^{19} + 7 T^{18} + 
5 T^{17} + 4 T^{16} + 2 T^{15} - 
T^{14} - 3 T^{13} - 5 T^{12} - 5 T^{11} \\
- 5 T^{10} - 5 T^{9} - 3 T^8 - T^7
+ 2 T^6 + 4 T^5 + 5 T^4
+ 7T^3 + 6 T^2 + 5 T  + 3).
\end{gather*}
Not only are the two assertions consistent, but the characteristic
polynomial computed from $p$-adic cohomology demonstrates the symmetry
forced by the functional equation of the zeta function, i.e., by Poincar\'e
duality on cohomology. (The geometric Picard number in this instance 
appears to be 4, as the characteristic polynomial appears to be divisible
not only by $T+1$ but also by $T^2+1$, so Tate's conjecture
predicts extra cycle classes
over $\FF_{3^2}$ and again over $\FF_{3^4}$.)
\end{remark}

\begin{remark}
It seems an interesting question to explore to what extent one can recover
a zeta function from the various pieces of data we have in the above
situation:
\begin{itemize}
\item
the symmetry and location of roots, from the Weil conjectures;
\item
the initial point counts;
\item
divisibilities implied by the relationship between the Newton and 
Hodge polygons
(see Remark~\ref{R:char poly});
\item
congruences derived from computing $p$-adic cohomology to low precision.
\end{itemize}
In particular, in many cases it may be possible to combine information
to recover zeta functions using $p$-adic cohomology calculations at much
less precision than would be predicted by a straightforward application
of the Weil conjectures plus taking into account the Hodge polygon.
In the case of Example~\ref{exa:f3 quartic}, this is discussed in
detail in \cite{kedlaya-search}; the end result is that 
not only is the zeta function equal to the guess
of Remark~\ref{R:compare charpoly}, but in fact this already follows
from the computation of the characteristic polynomial of $3^{-1} F$
modulo $3^1$ (and so from the Frobenius matrix with final precision
$3^3$) without counting any rational points at all.
\end{remark}

\subsection{Examples: degree 4 over $\FF_p$ ($p = 5,\dots,19$)}

We next exhibit some examples where we can compute the geometric Picard
number. These examples appear in a construction of van Luijk
\cite[Proposition~5.1]{vanluijk}; for context, 
we first state  \cite[Proposition~4.1]{vanluijk}.
\begin{prop}[van Luijk] \label{P:van luijk}
Let $k$ be a field, and suppose $\alpha,\beta \in k$ satisfy
$\alpha^3 \beta \neq \alpha \beta^3$. Let $f \in k[w,x,y,z]$
be a homogeneous polynomial of degree $3$, such that either the
coefficients of $y^3$ and $z^3$ in $f$ differ, or the coefficients of
$x^2 y$ and $x^2 z$ in $f$ differ. Suppose that the surface $X$ in
$\PP^3_k$ defined by 
\begin{equation} \label{eq:vanluijk}
wf - (xy + xz + \alpha yz)(xy + xz + \beta yz)
\end{equation}
is smooth with geometric Picard number $2$, and put
$\overline{X} = X \times_k \overline{k}$. Then
the group $\Aut(\overline{X})$ is trivial.
\end{prop}
Note that the surface in \eqref{eq:vanluijk}, if smooth, has arithmetic
Picard number at
least 2, as the hyperplane section $w=0$ splits into two conics.

\begin{example}
Over $\FF_5$, we verify the instance of Proposition~\ref{P:van luijk}
with $\alpha = 1$, $\beta = 3$, and
\[
f = 3x^3 + 3xy^2 - xyw + 3xzw - xw^2 + y^3 - y^2w + 2z^3 + w^3;
\]
one case of 
\cite[Proposition~5.1]{vanluijk} relies on the fact that this surface
has geometric Picard number at most 2 (and hence has Picard number and 
geometric Picard number exactly 2). This we verify in turn by computing
the matrix modulo $5^3$, which requires initial precision $5^7$,
5482 CPU seconds, and 514 MB of memory; applying
Algorithm~\ref{algo:picard} to bound the contributions to the right side
of \eqref{eq:tate4} shows that the geometric Picard number is indeed at most 2.
\end{example}

For $p>5$, there are relatively few values of $n$ that can contribute
to \eqref{eq:tate4} for which $\zeta_n$ is $p$-adically close to 1, so 
one might expect that less final precision is needed to bound the geometric
Picard number. This expectation is confirmed by the following set of examples.

\begin{example} \label{exa:larger p}
Table~\ref{tab:deg4fp} lists 
some further instances in which the geometric Picard number condition
in Proposition~\ref{P:van luijk} can be verified
(for $p=7, 11, 13, 17, 19$), yielding the remaining cases of
\cite[Proposition~5.1]{vanluijk}. We also list two cases
where the geometric Picard number condition appears to hold
but has not been verified definitively
(for $p=23, 29$). In all cases, we took
$\alpha = 1$, $\beta = 3$, and final precision $p^2$; the initial precision
required was always $p^4$, and for $p \leq 19$,
the upper bound in \eqref{eq:tate4}
was found to be 2. The examples for $p=11, 13,17,19$
were found by prescreening with initial precision $p^3$;
the examples for $p=23, 29$ were obtained by prescreening, but 
we did not complete the calculation at initial precision $p^4$.
\begin{table}[ht]
\caption{Further instances and presumed instances of Proposition~\ref{P:van luijk} over $\FF_p$.}
 \label{tab:deg4fp}
\begin{center}
\begin{tabular}{|c|c|c|c|}
\hline
$p$ & $f$ & CPU sec & MB \\
\hline
7 & 
$-2 x^3 + 2x^2y + 2x^2w + y^3 + 3y^2w + yzw - yw^2 + 2z^3$ &
697 & 63 \\
11 & $-5x^3 - 2x^2y -5 xy^2 - 2xz^2 + y^3 - yz^2 + 2z^3 - 4w^3$
& 5320 & 106 \\
13 & $3x^3 -6 x^2 z + y^3 -6 y w^2 + 2 z^3$ & 14997 & 158 \\
17 & $- x^3 + 8 x^2 y - x y w + y^3 - y^2 w + 2 z^3 + 5 z^2 w 
    -5 z w^2$ & 61996 & 306 \\
19 & 
$6 x^3 + 3 x^2 z + 7 x y w -7 x z^2 + 8 x z w -9 x w^2 + y^3$ & 116323 & 459 \\
& $- y^2 w -5 y z^2 + 5 y w^2 + 2 z^3 -4 z^2 w -2 z w^2$ &  &  \\
\hline \hline
23 & $-11 x^3 -9 x^2 y -2 x^2 z -5 x^2 w -3 x y z$ & ? & ? \\
&  $-10 x y w + y^3 +     11 y w^2 + 2 z^3 -4 w^3$ &  &  \\
29 & $4 x^3 + 4 x^2 w -5 x y^2 -14 x y w + y^3$ & ? & ? \\
&  $+ 7 y^2 z -3 y z^2 + 3 y w^2 + 2 z^3 - w^3$ &  &  \\
\hline
\end{tabular}
\end{center}
\end{table}
\end{example}

\subsection{Examples: degrees 4, 5 over $\FF_2$}
\label{subsec:f2}

To conclude, we edge towards the realm of coding theory proper,
 by considering examples over $\FF_2$.
Frustratingly, the asymptotic advantage obtained in the $p$-adic
algorithms by taking $p = 2$ (which occur because the degree of the
truncated Frobenius lift is a linear function of $p$) is somewhat 
counterbalanced by the need for additional precision to fight the
unfortunate propensity of small integer denominators to be divisible by large
powers of 2. Nonetheless, one can still say something.

\begin{example}
The smooth quartic surface over $\FF_2$ defined by the polynomial
\[
x^4 + x^3 z + x^2 y^2 + x^2 y w + x^2 z^2 + x^2 z w + 
    x y^3 + y^4 + y^3 w + y z^3 + z^4 + z^2 w^2 + w^4
\]
was found to have Picard number 1 by computing a Frobenius matrix
modulo $2^4$; this required initial precision $2^{13}$, 7182 CPU seconds,
and 472 MB of memory. 
This example was found by prescreening using initial precision
$2^8$, which required a mere 88 CPU seconds and 52 MB of memory.
\end{example}

If we pass from quartics to quintics, then middle cohomology
becomes 53-dimensional, with Hodge numbers 4, 45, 4. (Again, the space
we compute in is only 52-dimensional because passing from $\tilde{X}$ to
$\tilde{U}$ removes one cycle class.) This removes the parity
obstruction to having \emph{geometric} Picard number 1, and one would expect
to be able to find examples of such. Moreover, since in this case one is not
forced to include the eigenvalue $-2$, which is indistinguishable from $+2$
modulo $2^2$, one might even expect to be able to work with less precision.
The following example fulfills these expectations.

\begin{example} \label{exa:f2 quintic}
The smooth quintic surface over $\FF_2$ defined by the polynomial
\begin{gather*}
 x^5 + x^3 y z + x^2 y^2 w + x^2 y z^2 + x^2 z^3 +
    x z^2 w^2 \\ +y^5 + y^3 z w + y^2 z w^2 + y z w^3
    + z^5 + z^2 w^3 + w^5
\end{gather*}
was found to have geometric Picard number 1 
by computing a Frobenius matrix modulo
$2^3$; this required initial precision $2^{12}$,
22685 CPU seconds, and 179 MB of memory.
Note that in this case screening for the geometric Picard number is a
nontrivial calculation, because we must check up to $n=210$, which entails
working in some large extensions of $\QQ_2$.

The example was found by prescreening with initial precision $2^6$.
If we had needed final precision $2^4$, we would have required
initial precision $2^{13}$; we project that in our implementation,
such a calculation would
require roughly $100000$ CPU seconds, and would have to be done in
batches (or parallelized) to avoid memory overrun.
\end{example}

\section*{Acknowledgments}

Kedlaya wishes to thank
 Felipe Voloch, whose discussions about his work with Zarzar
initiated this project, Ronald van Luijk for helpful
discussions about \cite{vanluijk}, and Johan de Jong and Bas Edixhoven 
for additional discussions. 
Some of these discussions took place at the Oberwolfach workshop
``Explicit methods in number theory'' in July 2005, where some of this
material was presented in early form. 
Abbott and Roe were supported by
MIT's Undergraduate Research Opportunities Program.
Kedlaya was supported by NSF grant DMS-0400747.

\end{document}